\documentclass[10pt]{amsart}

\usepackage[utf8]{inputenc}
\usepackage[T1]{fontenc}
\usepackage{lmodern}

\usepackage{amssymb,appendix}
\usepackage{mathtools}
\usepackage{amscd}
\usepackage[mathscr]{eucal}
\usepackage{enumitem}
\usepackage[hmargin=1in,vmargin=1in]{geometry}

\usepackage{bk_macros}

\usepackage{graphicx}
\usepackage{xcolor}

\usepackage[roman]{sublabel}

\usepackage[labelfont=bf,font=small]{caption}
\usepackage{subcaption}

\usepackage[noadjust]{cite}
\usepackage[foot]{amsaddr}

\usepackage{placeins}
\usepackage{bm}
\usepackage{siunitx}

\definecolor{cite_purple}{RGB}{128,9,158}
\definecolor{cite_blue}{RGB}{2,95,176}
\definecolor{link_red}{rgb}{0.7,0,0}
\usepackage[colorlinks=true,
            citecolor=cite_purple,
            linkcolor=link_red,
            urlcolor=cite_blue]{hyperref}

\newlength{\arrayrulewidthOriginal}

\theoremstyle{plain}

\theoremstyle{definition}

\begin{document}

\title{Dynamic Homeostasis in Relaxation and Bursting Oscillations}

\author{Christopher J. Ryzowicz$^{1}$}
\address{$^{1}$ Department of Mathematics, Florida State University, Tallahassee, Florida, 32306, USA}

\author{Richard Bertram$^{1,2,3}$}
\address{$^{2}$Institute of Molecular Biophysics, Florida State University, Tallahassee, Florida 32306, USA}
\address{$^3$ Program in Neuroscience, Florida State University, Tallahassee, Florida 32306, USA}

\author{Bhargav R. Karamched$^{*,1,2,3}$}
\thanks{$^{*}$\href{mailto:bkaramched@fsu.edu}{bkaramched@fsu.edu} (Corresponding author)}



\date{\today}

\begin{abstract}
\textbf{Relevance to Life Sciences}\\
Homeostasis, broadly speaking, refers to the maintenance of a stable internal state when faced with external stimuli. Failure to manage these regulatory processes can lead to different diseases or death. Most physiologists and cell biologists around the world agree that homeostasis is a fundamental tenet of their disciplines. Nevertheless, a precise definition of homeostasis is hard to come by. Often times, homeostasis is simply defined as ``you know it when you see it''. 
\newline\noindent\textbf{Mathematical Content}\\
Mathematical treatments of homeostasis involve studying equilibria of dynamical systems that are relatively invariant with respect to parameters. However, physiological processes are rarely static and often involve dynamic processes such as oscillations. In such dynamic environments, quantities such as average values may be relatively invariant with respect to parameters. This has been referred to as ``homeodynamics''. We present a general framework for homeodynamics involving systems with two or more time scales that elicits homeostasis in the temporal average of a species. The key point is that homeostasis manifests when measuring the slow variable responsible for driving oscillations and is not apparent in the fast variables. We demonstrate this in the Fitzhugh-Nagumo model for relaxation oscillations and then in two models for electrical bursting activity and calcium oscillations in pancreatic $\beta$-cells. One of these models has multiple slow variables, each driving the bursting oscillations in different parameter regimes, but homeodynamics is only present in the variable currently engaged in this role. \\
\noindent\textbf{Further Reading}\\
\noindent \textit{Homeostasis} \cite{guyton2006text,reed2017analysis,nijhout2014escape,duncan2018homeostasis,nijhout2019systems}\\
\noindent \textit{Dynamic homeostasis} \cite{lloyd2001homeodynamics,xiong2023physiological,yates1994order}\\
\noindent \textit{Relaxation oscillations} \cite{bertram2017multi,hastings18}\\
\noindent \textit{Bursting oscillations} \cite{rinzel1987formal, bertram1995topological,rinzel1986different,rinzel1987dissection,terman1992transition,bertram2000phantom,bertram2004calcium,bertram2023deconstructing,desroches22,rinzel98}

\end{abstract}

\maketitle

\raggedbottom
\thispagestyle{empty}

\section{Introduction}
A central question in physiology and cell biology is how various systems maintain homeostasis. Homeostasis is identified as one of the key tenets of physiology and biology by many top researchers around the world~\cite{michael2007conceptual,alpern2011competencies}. Surprisingly, its precise definition and the biomolecular control systems that facilitate it remain an open question~\cite{modell2015physiologist}. 

Why is homeostasis important? Broadly speaking, homeostasis refers to the maintenance of a stable internal state even when faced with (potentially drastic) external disturbances. Homeostatic conditions in humans typically include a core body temperature of 98.6$^{\circ}$F (37$^{\circ}$C), blood glucose levels around 85 mg/dL (4.7 mmol/L), and a slightly basic blood pH of about 7.4.
Failure to regulate processes so that such values are maintained in a multitude of environments can lead to severe disorders or death~\cite{guyton2006text}. 

\begin{figure}[h!]
\centering
\includegraphics[width=8 cm, height= 8 cm]{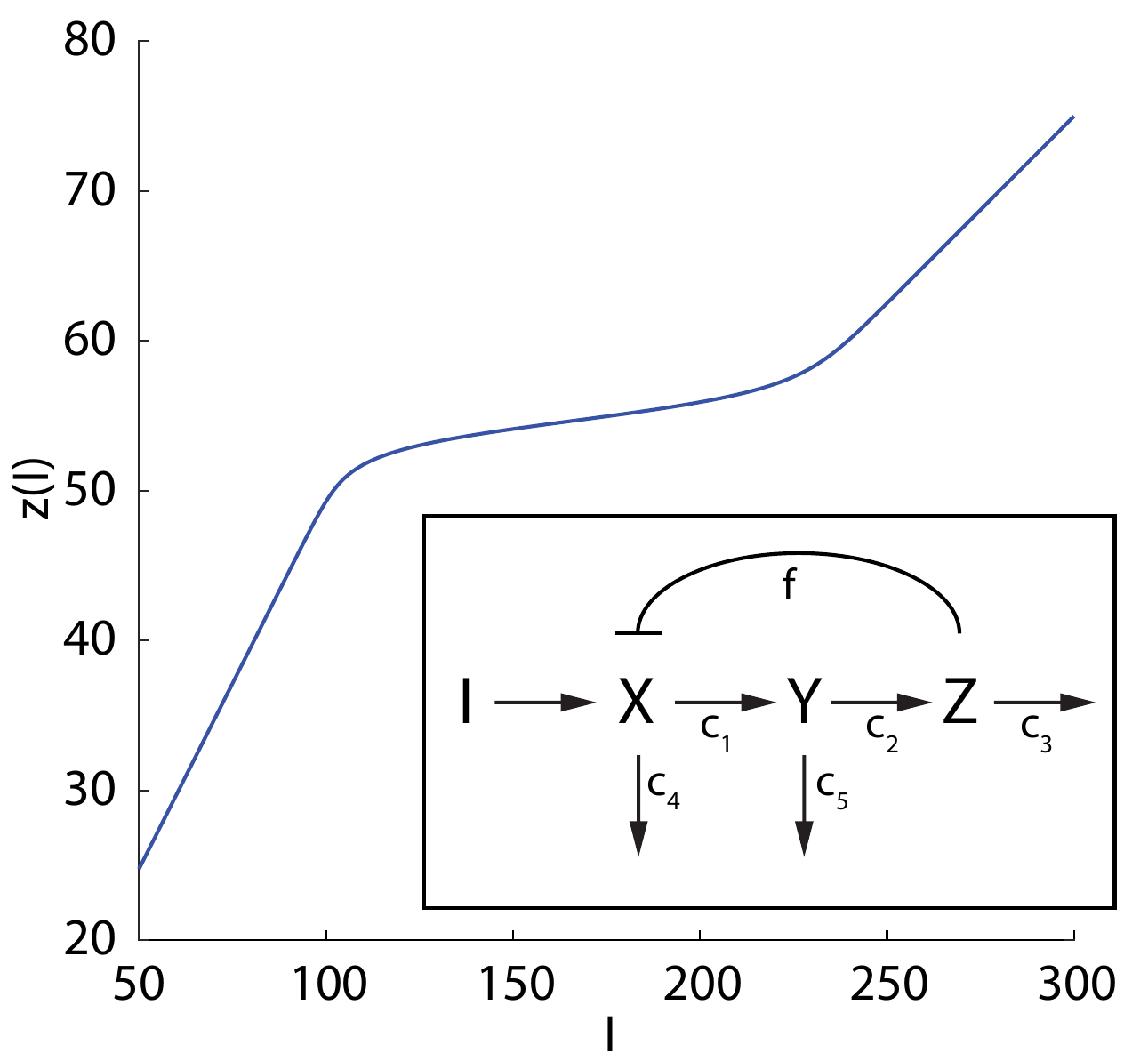}
\caption{Homeostatic chair curve of equilibria $z$ as a function of the input values $I$ recreated from \cite{reed2017analysis} using the feedback inhibition model (inset).}
\label{fgr1: Reed Chair}
\end{figure}

To understand homeostasis, physiologists and cell biologists have adopted the theory and vernacular of control systems---thereby framing complicated biophysical networks with a mathematical lens. Mathematically, homeostasis can be interpreted as a circumstance wherein the equilibrium of a dynamical system is relatively invariant with respect to alterations to the system parameters. Indeed, Reed and co-authors published a series of papers that investigate various biomolecular control systems to determine what network motifs facilitate homeostasis~\cite{reed2017analysis,nijhout2014escape,nijhout2019systems,golubitsky2017homeostasis}. An example of a network motif facilitating homeostasis is feedback inhibition (see Fig.~\ref{fgr1: Reed Chair} (inset)). In this example, the input, $I$, feeds into a pathway that amplifies biochemical species $X$. This leads to the eventual production of biochemical species $Z$, which inhibits the action of $X$. In Fig.~\ref{fgr1: Reed Chair}, we show the equilibrium value of $Z$ as $I$ is varied. The curve resembles a chair, with the ``seat'' of the chair occurring over the interval of $I$ values for which the equilibrium $Z$ values are relatively invariant, which we refer to as semi-invariant. This semi-invariance is representative of homeostasis.

Golubitsky, Stewart, and others have developed rigorous mathematical theory characterizing the conditions under which a system exhibits homeostasis~\cite{golubitsky2017homeostasis,golubitsky2018homeostasis,golubitsky2020input,wang2021structure,duncan2019coincidence,duncan2018homeostasis}. The basic idea is that homeostasis occurs when the derivative of an output with respect to the input is zero. To precisely characterize homeostasis, they invoke singularity theory.

A limitation of these characterizations of homeostasis is the focus on  equilibria. Many physiological observables exhibit rich dynamical profiles such as oscillations, bursting~\cite{rinzel1987formal,bertram1995topological,izhikevich2000neural}, or even chaos~\cite{karamched2021delay}. Quantities commonly associated with static equilibrium values, such as body temperature or blood glucose levels, actually oscillate \cite{aschoff1971human,sturis1991computer,xiong2023physiological}. This suggests that more complicated temporal signatures, not just static equilibria, may encode homeostasis.


In this paper, we propose a general framework for a novel physiological architecture of homeostasis, called dynamic homeostasis or homeodynamics \cite{lloyd2001homeodynamics,xiong2023physiological}. We focus on excitable systems with multiple timescales that undergo bifurcations to limit cycles. For some input values, such systems will rapidly contract to a stable equilibrium. However, the equilibrium can destabilize via, for example, a Hopf bifurcation or homoclinic bifurcation. With multiple timescales, the resulting oscillation will manifest as a relaxation oscillation characterized by a slow build-up process followed by fast transitions \cite{bertram2017multi}. Such relaxation oscillations are observed in biological rhythms of calcium \cite{sneyd2017dynamical}, voltage \cite{morris1981voltage}, and cardiac dynamics \cite{franzone2006adaptivity}.

We show that homeostasis is observed in relaxation oscillations when looking at the temporal average of the oscillatory signals. In particular, the temporal average of the slow variable exhibits homeostasis---the fast variable fails to exhibit homeostasis. The underlying reason for this is that the trajectory of the limit cycle in phase space is semi-invariant with respect to input parameters. Specifically, the slow variable oscillates back and forth as a sawtooth wave between the same values for all oscillation-facilitating input values. Thus, its temporal average does not change. The fast variable appears as a square wave that jumps between the same values for all inputs in the oscillatory regime; however, the duty cycle of the time signal varies, thereby yielding different temporal averages. We thus show how a complicated temporal signal may encode homeostasis and establish which variables to measure to detect it.

We first demonstrate this phenomenon in the Fitzhugh-Nagumo model. We further show that including noise in the input allows for the relaxation oscillator to work even better as an encoder of homeostasis. We then generalize our findings to a biologically relevant bursting oscillation using both a deterministic and stochastic version of the Chay-Keizer model for pancreatic $\beta$-cells \cite{chay1983minimal}. Thereafter, we briefly look at a system with multiple slow timescales \cite{bertram2004calcium} to show how homeostasis may be observed in systems with multiple slow variables.



\section{Relaxation Oscillations with the FitzHugh-Nagumo model}

The FitzHugh-Nagumo (FHN) model is a simple, two-dimensional model, adapted from the van der Pol oscillator \cite{vanDerPol1927ForcedOscillations}. Introduced by Richard FitzHugh \cite{fitzhugh1961impulses} in 1961 and refined by Jinichi Nagumo \cite{nagumo1962active} a year later, the model is also a simplified version of the 4-dimensional Hodgkin-Huxley equations for neuronal spiking \cite{HodgkinHuxley1952quantitative}. It replicates many neuronal phenomena from the Hodgkin-Huxley equations but is easier to analyze due to its lower dimensionality.
\subsection{Excitability and relaxation oscillations in the FHN model}

\begin{figure}[h!]
\centering
\includegraphics[width= \columnwidth]{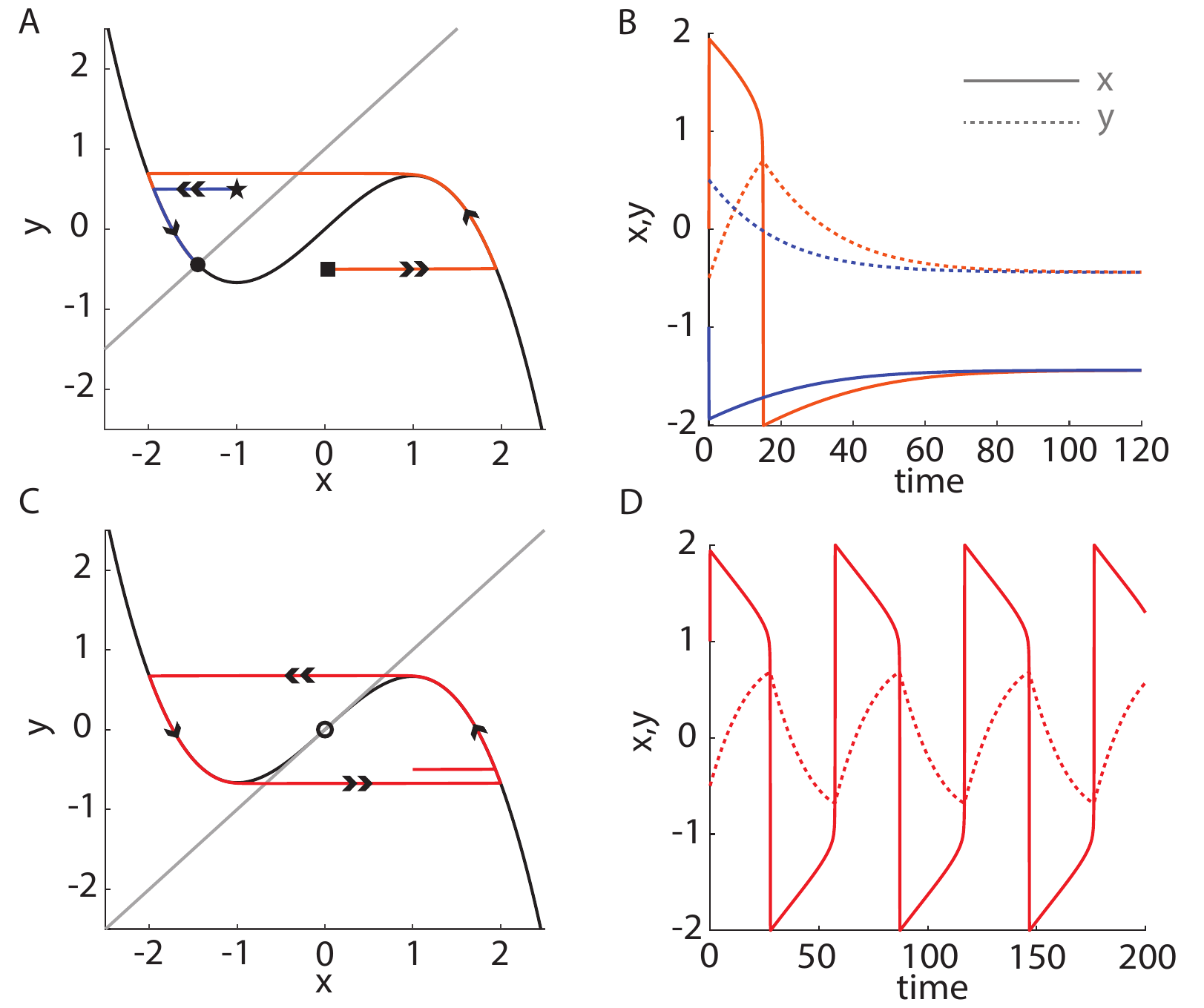}
\caption{Dynamics of system (\ref{eq1: FHN Model}) for $\alpha=1$, $\mu=30$ and (A-B) $J=1$ or (C-D) $J=0$. Left column: dynamics in the $(x,y)$-plane with $x$-nullcline (black), $y$-nullcline (gray) and sample trajectories (color). The blue curves indicate sub-threshold responses to perturbation from rest, while the orange curves indicate super-threshold responses. The double arrows indicate fast motion, while the single arrows indicate slower motion. The filled circle is a stable equilibrium, the open circle is an unstable equilibrium. Right column: sample time series showing (B) sub-threshold and super-threshold responses, and (D) relaxation oscillations.}
\label{fgr2: FHN Excitability and Oscillations}
\end{figure}

We begin by reviewing two fundamental features of the FHN model: excitability and oscillations. The  nondimensional model satisfies the ordinary differential equations
\begin{equation}\label{eq1: FHN Model}
    \frac{dx}{dt} = \mu(x-\frac{x^3}{3}-y) \qquad \frac{dy}{dt} = \frac{1}{\mu}(J+\alpha x-y),
\end{equation}
where $x(t)$ and $y(t)$ are the activator variable (analogous to the membrane potential) and the recovery variable (analogous to a K$^+$ channel gating variable), respectively \cite{fitzhugh1969mathematical}. The parameter $\alpha>0$ controls the strength of coupling feedback, and the parameter $\mu\gg1$ controls the timescale separation (we use $\mu =30$ in all simulations). We consider the case in which the recovery variable is modulated by an input parameter $J$ that could reflect gene expression of the ion channels responsible for spike repolarization.

The equilibria of the FHN model are obtained by setting both temporal derivatives in Eqs.~\eqref{eq1: FHN Model} equal to zero. This results in the following algebraic equations for the equilibria $(x^*,y^*)$:
\begin{equation*}
    \frac{{x^*}^3}{3} + (\alpha-1)x^*+J=0, \qquad y^*=J+\alpha x^*.
\end{equation*}
For $\alpha\geq1$ the system admits one real equilibrium. If the input parameter $J$ satisfies either $J<\frac{2}{3}-\alpha$ or $J>\alpha-\frac{2}{3}$, the equilibrium will be locally stable (see Appendix \ref{Appendix: FHN Linear stability and Hopf bifurcation} for linearization and stability analysis). Fig. \ref{fgr2: FHN Excitability and Oscillations}A-B depicts the possible dynamics of this system when the equilibrium is stable and $\mu$ is sufficiently large. The stable node (filled black circle) occurs at the intersection of the $x$-nullcline (black cubic curve) and the $y$-nullcline (gray line). If we impart a small perturbation away from the equilibrium in one direction (black star), a sub-threshold response occurs. That is, the system passively converges back to equilibrium. However, if we impart a perturbation in a different direction (black square), a super-threshold response occurs. The trajectory undergoes a `spike' before converging back to equilibrium. The threshold between these two different behaviors is the middle branch of the cubic $x$-nullcline. This threshold phenomenon, called excitability, is the defining feature of a model for impulse generation. 

\begin{figure}[h!]
\centering
\includegraphics[width=\columnwidth]{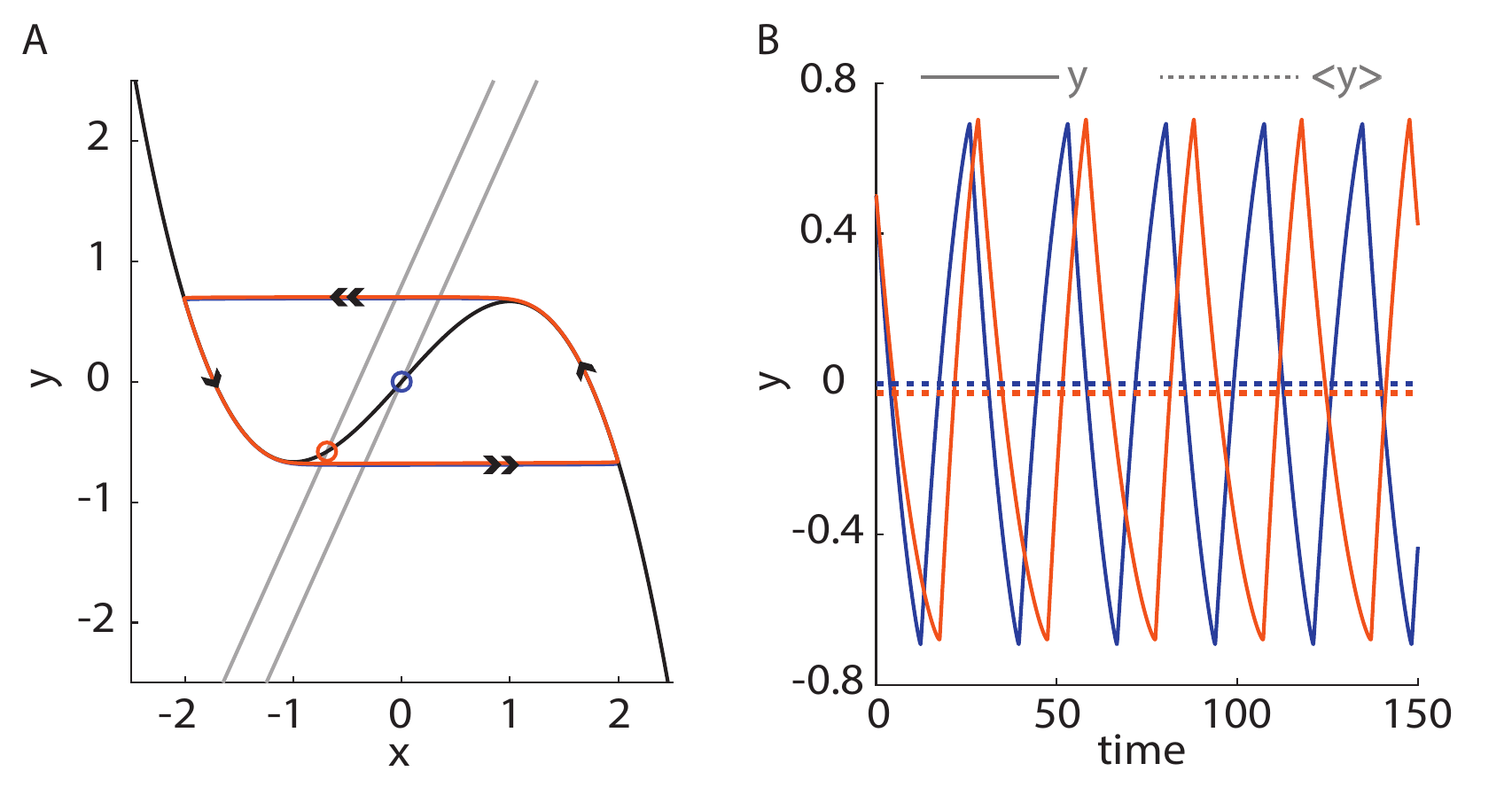}
\caption{Semi-invariance in the mean of the slow variable with changes in the input. (A) Phase portrait for the case $J=0$ (blue) and $J=0.8$ (orange). Increasing $J$ translates the $y$-nullcline (gray lines) upward and leaves the $x$-nullcline (black curve) unchanged. This results in a left shift in the equilibrium, but no change in stability (open circles). Limit cycles for the two cases (colored curves) are almost identical. (B) Time series for $y$ in the two cases (solid curves) and the time-average values (dashed lines) show semi-invariance of the mean value of the slow variable. For both cases, $\alpha=2$.}
\label{fgr3: FHN Invariance}
\end{figure}

When $J \approx \pm(\frac{2}{3}-\alpha)$, the system undergoes a subcritical Hopf bifurcation (see the Appendix \ref{Appendix: FHN Linear stability and Hopf bifurcation} for Hopf bifurcation analysis), rendering the equilibrium unstable for $\frac{2}{3}-\alpha < J < \alpha-\frac{2}{3}$ and causing the onset of periodic solutions (see Fig. \ref{fgr2: FHN Excitability and Oscillations}C). These emergent periodic solutions are unstable, but quickly change stability and grow in amplitude to become relaxation oscillations that are characterized by a substantial timescale separation between the state variable dynamics due to the choice of the parameter $\mu$ (Fig. \ref{fgr2: FHN Excitability and Oscillations}D). In this model, $x$ evolves on a fast timescale and $y$ evolves on a slower timescale. See \cite{bertram2017multi} for discussion of relaxation oscillations.

\subsection{Slow variable time-average shows semi-invariance, fast variable does not} 

For the remainder of the paper we focus primarily on persistent oscillations, starting with relaxation oscillations produced by the FHN model. Examples of two relaxation oscillations, with $J=0$ and $J=0.8$, are shown in Fig. \ref{fgr3: FHN Invariance}.
Increasing $J$ from 0 to 0.8 results in an upward translation of the $y$-nullcline, shifting the equilibrium closer to the lower knee of the cubic $x$-nullcline (Fig. \ref{fgr3: FHN Invariance}A). This change has almost no effect on the limit cycle; the limit cycles for the two values of $J$ lie almost directly on top of each other in phase space. The time series for $y$ has a sawtooth shape for both values of $J$ (Fig. \ref{fgr3: FHN Invariance}B). Although the limit cycle is almost unchanged by the change in $J$, the period is somewhat different, so the $y$ oscillations drift apart over time. In spite of this, the time-average of $y$, $\langle y \rangle$, is almost identical for the two values of the input. 

\begin{figure}[h!]
\centering
\includegraphics[width=\columnwidth]{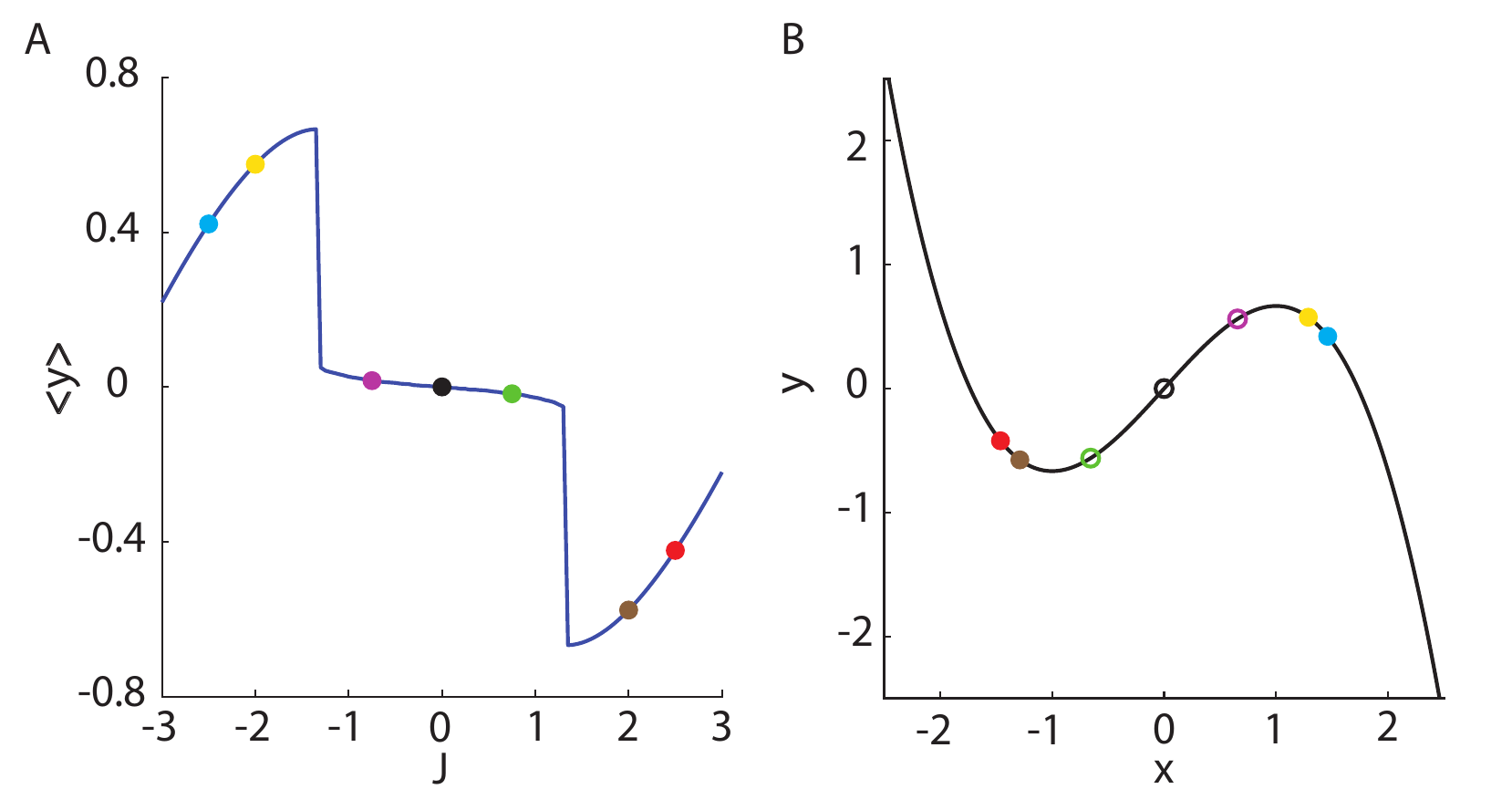}
\caption{Dynamic homeostasis in the time-average of the slow variable of the FHN model. (A) Homeostatic chair over a range of input values $J \in [-3,3]$ with $\alpha=2$. (B) Phase space showing the $x$-nullcline (black curve) and colored points indicating the locations of equilibria with different $J$ values. The equilibria are color coded to correspond with the colored points in panel A.}
\label{fgr4: FHN Chair}
\end{figure}

To determine the extent of the semi-invariance, we simulated the FHN model and computed $\langle y \rangle$ over a range $J \in [-3,3]$ increasing by increments of $0.05$. The result is shown in Fig. \ref{fgr4: FHN Chair}A. For small values of $J$, below the threshold for oscillations, increasing the input causes a substantial increase in $\langle y \rangle$. For large values of $J$, above the threshold for oscillations, increasing the input also causes a substantial increase in $\langle y \rangle$. For intermediate values of $J$, however, changes in the input lead to very little change in $\langle y \rangle$. This semi-invariance in $\langle y \rangle$ is an example of dynamic homeostasis \cite{xiong2023physiological,yates1994order,soodak1978homeokinetics,lloyd2001homeodynamics}.
This phenomenon results in a ``homeostatic chair'' curve similar to Fig.~\ref{fgr1: Reed Chair} and described in prior studies \cite{reed2017analysis,nijhout2014escape,duncan2018homeostasis,golubitsky2017homeostasis,nijhout2019systems}. 
The key difference, however, is that in those models the homeostasis is in the value of a variable at equilibrium. In our case, it is in the average value of a variable over an oscillation period, or dynamic homeostasis. As shown in Fig. \ref{fgr4: FHN Chair}, semi-invariance in the variables at equilibrium does not occur in the FHN model. 

The dynamic homeostatic plateau section (the ``seat'') of the chair curve exists in the same parameter range of $J$ in which the relaxation oscillations occur. Fig. \ref{fgr4: FHN Chair}B shows phase space with the $x$-nullcline (black curve) and seven points that indicate different equilibria for distinct values of $J$. The corresponding $\langle y \rangle$ values are shown in Fig.~\ref{fgr4: FHN Chair}A in matching colors. The three points on the seat (magenta, black, and green) correspond to the points which lie in between the two knees of the $x$-nullcline. These are locations where the equilibrium is unstable and oscillations occur as seen in Fig. \ref{fgr3: FHN Invariance}A. Outside of this range, the equilibria are locally stable (red, brown, yellow, and cyan) and the system does not oscillate. 

This observation, that the semi-invariant portion of the chair corresponds to the parameter interval in which there are relaxation oscillations, also explains the two large jumps in the chair curve. At a small value of $J$ there is a stable equilibrium with a large $y$ value (either the yellow or cyan point in Fig. \ref{fgr4: FHN Chair}). For a somewhat larger $J$ value for which the stable equilibrium is replaced by a stable periodic solution (magenta point), $y$ reaches down to negative values during the oscillation, yielding a $\langle y \rangle$ value much smaller than when the equilibrium was stable. Thus, there is a sudden jump down in the left portion of the chair curve, at $J \approx -1.33$. At the other extreme, when $J$ is positive and just past 1, the transition from periodic solutions (green point) covering both positive and negative $y$ values to a stable equilibrium with a negative $y$ value (brown point), results in a large decrease in $\langle y \rangle$. This is the second jump down, in the right portion of the chair curve at $J \approx 1.33$.


\subsubsection{There is no invariance in the time-average of the fast variable}
\begin{figure}[h!]
\centering
\includegraphics[width=\columnwidth]{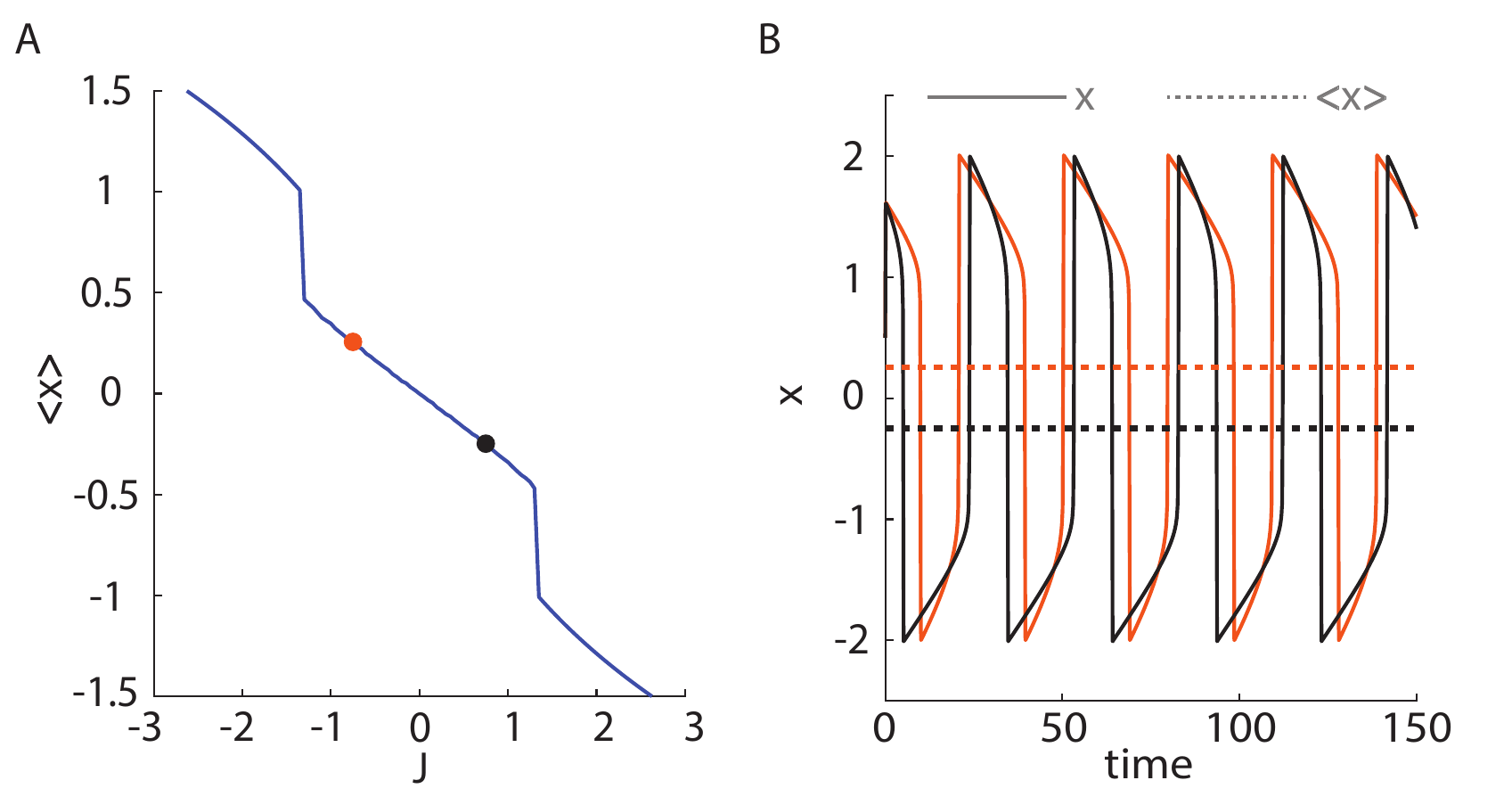}
\caption{Lack of dynamic homeostasis in the fast variable dynamics of the FHN model. (A) The time-average of $x$ plotted as a function of input $J$ shows no homeostasis. (B) The time series for two values of $J$ (solid curves) and time-averages (dashed lines) show a substantial difference in $\langle x \rangle$ due to the difference in the oscillation duty cycle. For both panels, $\alpha=2$.}
\label{fgr5: FHN loss of invariance in fast variable}
\end{figure}

We next consider how the time-average of the fast variable $x$, $\langle x \rangle$, changes with changes in $J$. This is plotted as a function of $J$ in Fig. \ref{fgr5: FHN loss of invariance in fast variable}A. As with the slow variable, there is no invariance in this variable over intervals of $J$ where the equilibrium is stable. Unlike before, however, there is no invariance in $\langle x \rangle$ over intermediate values of the input (approximately $-1.33<J<1.33$). The reason for this is evident in the $x$ time series (Fig. \ref{fgr5: FHN loss of invariance in fast variable}B), using two values of $J$ within the oscillation region (orange and black points in Fig. \ref{fgr5: FHN loss of invariance in fast variable}A). 
With the lower input value, the duty cycle (the duration of an ``up phase'' versus the entire oscillation period) is higher than with the higher input value. As a result, the time-average value of $x$ is higher with the lower input (dashed orange line) than with the higher input (dashed black line). The change in the duty cycle, which is responsible for the lack of invariance in $\langle x \rangle$, occurs because $J$ changes the proximity of the $y$-nullcline (gray curve in Fig. \ref{fgr3: FHN Invariance}A) to the left and right branches of the $x$-nullcline (black curve in Fig. \ref{fgr3: FHN Invariance}A). At a higher value of $J$ the $y$-nullcline is closer to the left branch of the $x$-nullcline than it is when $J$ is lower, so the speed of the trajectory along that branch is slower. As a result, when $J$ is larger there will be a longer ``down state'' (and a shorter ``up state'') in the $x$ oscillations, yielding a lower duty cycle. The central point here is that biological relaxation oscillations can exhibit dynamic homeostasis in the mean, but only in the slow variable.


\subsection{The homeostatic range depends on the feedback strength}
\begin{figure}[h!]
\centering
\includegraphics[width=\columnwidth]{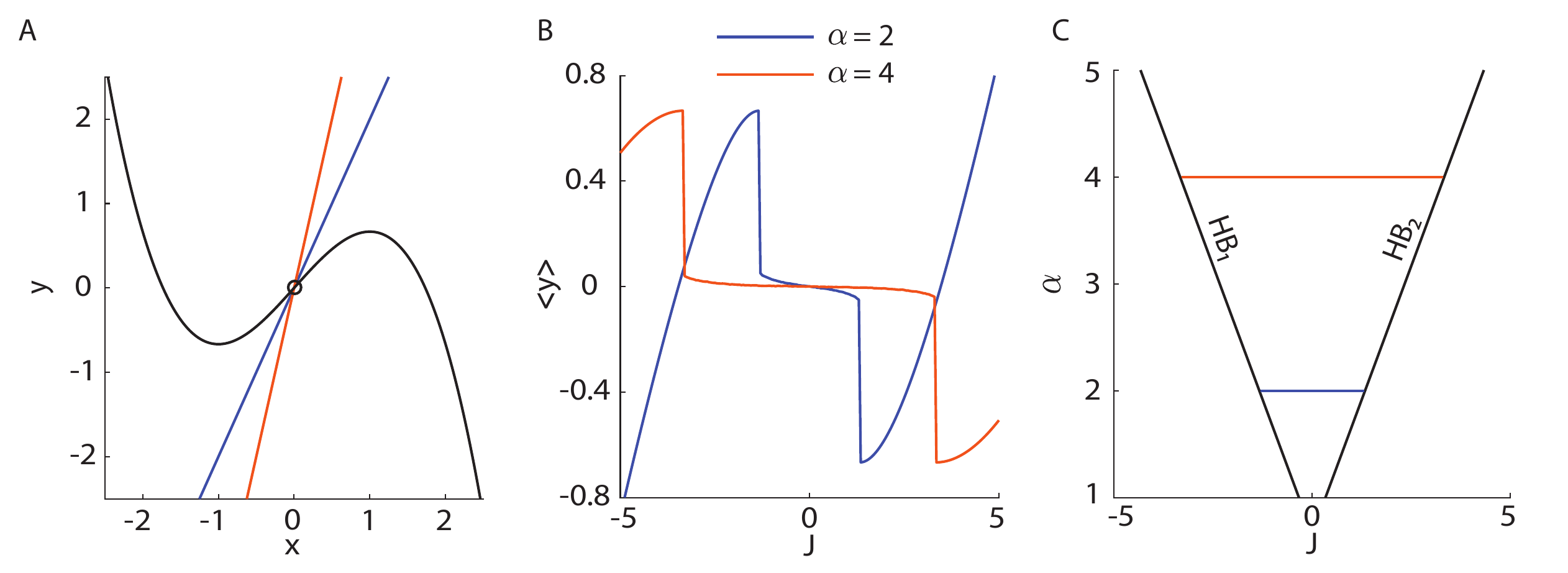}
\caption{Semi-invariance depends on the feedback strength. (A) Phase plane when $\alpha=2$ (blue) and $\alpha=4$ (orange). The y-nullcline slope increases with $\alpha$ while x-nullcline (black curve) remains unchanged. (B) Homeostatic chair curves of $\langle y \rangle$ as a function of $J$ when $\alpha=2$ (blue) and $\alpha=4$ (orange). (C) Two parameter bifurcation diagram in $(J,\alpha)$-space with Hopf bifurcation solutions (black curves). }
\label{fgr6: FHN invariance dependance on alpha}
\end{figure}

The feedback strength, $\alpha$, affects the slope of the $y$-nullcline, but has no effect on the $x$-nullcline of the FHN model. In particular, increasing $\alpha$ increases the slope of the $y$-nullcline (Fig. \ref{fgr6: FHN invariance dependance on alpha}A). What effect does this have on the interval of dynamic homeostasis? To answer this question we simulated the FHN model over a wide range of $J$ values and determined the asymptotic time-average of the slow variable for different $\alpha$ values. When $\alpha$ is increased from 2 (blue) to 4 (orange), the range of $J$ values at which there is dynamic homeostasis increases 
 (Fig. \ref{fgr6: FHN invariance dependance on alpha}B).  

The reason for the lengthening of the homeostatic interval can be explained by the effect of $\alpha$ on the location of the Hopf bifurcations that initiate and terminate the relaxation oscillations. In Fig. \ref{fgr6: FHN invariance dependance on alpha}C we plot a two parameter bifurcation diagram of the Hopf bifurcations in the $(J,\alpha)$-space (see Appendix \ref{Appendix: FHN Linear stability and Hopf bifurcation} for Hopf bifurcation analysis). The branches of the loci of Hopf bifurcation points (black) are approximately linear. The range of $J$ values in which there is a homeostatic plateau in  $\langle y \rangle$ is approximately the range of values of $J$ between the Hopf bifurcation branches. 
Thus, we conclude that the homeostatic region grows approximately linearly as $\alpha$ increases.

\subsection{Input stochasticity extends the range of homeodynamics in the FitzHugh-Nagumo model}
Biological systems, especially at the intracellular level, operate in crowded, heterogeneous, noisy environments \cite{tsimring2014noise,bressloff2014stochastic,losick2008stochasticity, kaern2005stochasticity}. Hence, it is important to develop biological models that incorporate stochasticity and analyze its effects on the dynamics. In this section, we analyze the effects of stochasticity in input on the homeostasis property of the FHN model.

The modified differential equations are 
\begin{equation}\label{eq2: FHN Stochastic Model}
    \frac{dx}{dt} = \mu(x-\frac{x^3}{3}-y) \qquad \frac{dy}{dt} = \frac{1}{\mu}(J^*+\alpha x-y),
\end{equation}
where all parameters and variables are defined as in (\ref{eq1: FHN Model}), except the input stimulus $J^*$ is a normally distributed random variable
\begin{equation}\label{eq3: FHN Normally Distributed J}
    J^* \sim N(J,\sigma),
\end{equation}
where $J$ is the mean of the normal distribution and $\sigma$ is the standard deviation (strength of the noise). 

\begin{figure}[h!]
\centering
\includegraphics[width=\columnwidth]{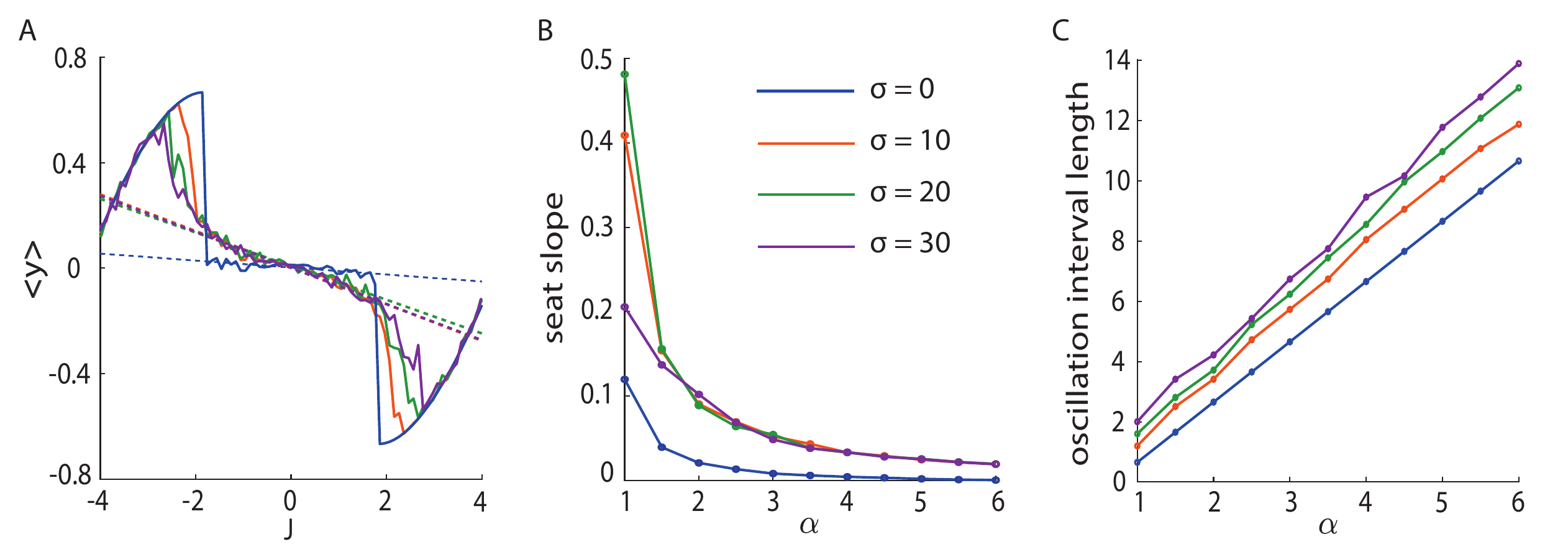}
\caption{Dynamic homeostasis is preserved when stochasticity is added to the input stimulus of the FHN model. (A) Homeostatic chair curve for the deterministic model ($\sigma=0$) and stochastic models ($\sigma=10,20,30$). Linear regression (dashed lines) show trends within the seat portion of the chair. A feedback strength of $\alpha=2.5$ was used for the simulations. (B) The seat slope declines with $\alpha$ in both the deterministic and stochastic models. (C) The range of $J$ values for which the FHN model exhibits oscillations increases with $\alpha$ for both the deterministic and the stochastic models.}
\label{fgr7: Stochastic FHN model}
\end{figure} 

We numerically solved (\ref{eq2: FHN Stochastic Model}) using the forward Euler method, randomly assigning the input stimulus every integral time step according to (\ref{eq3: FHN Normally Distributed J}). We then calculated the time-average of the slow variable, $\langle y \rangle$, after sufficiently long time for many different values of the normal distribution mean, $J$. The time-average as a function of $J$ for different standard deviations ($\sigma=0,10,20,30$) is given in Fig. \ref{fgr7: Stochastic FHN model}A. For very small and very large values of $J$, increasing (decreasing) the mean input results in an increase (decrease) in $\langle y \rangle$ regardless of the standard deviation. However for intermediate values of $J$, the effects of changing $J$ on $\langle y \rangle$ depends on the standard deviation of the distribution. We quantified these changes in two different ways: (1) measuring the slopes of the seats of the chair curves and (2) measuring the length of the oscillation interval in terms of $J$.

Dynamic homeostasis can be characterized as a segment of the $\langle y \rangle$ vs. $J$ curve that is flat, or has a slope near 0. This segment is like the seat of the chair, so we refer to it as the ``seat slope''. 
To determine the effects of noise in the input, we computed the seat slope for different values of the noise strength $\sigma$ using linear regression on the $\langle y \rangle$ values over the interval of $J$ between the two Hopf bifurcation values of the deterministic model (i.e. when $\sigma=0$). This was repeated for different values of $\alpha$, the parameter that sets the slope of the $y$-nullcline. The seat slope as a function of $\alpha$ is plotted in Fig. \ref{fgr7: Stochastic FHN model}B. One observation  is that the slope for the deterministic curve is less than that of all of the stochastic curves regardless of the $\alpha$ value. This is because the noise skews the chair curves near the ends of the seat (e.g. Fig. \ref{fgr7: Stochastic FHN model}A near $J \approx \pm 2$). It is also evident that the seat slopes of the stochastic model become very similar for all values of $\sigma$ as $\alpha$ increases. For low $\alpha$ values, the difference in slope values between the stochastic and deterministic curves is large because the distances between the Hopf bifurcations are small (see Fig. \ref{fgr6: FHN invariance dependance on alpha}B-C). That is, the seat intervals are short. For larger  $\alpha$ values, the seat slopes of the stochastic simulations not only converge to one another, but also approach those of the deterministic model. These numerical results show that dynamic homeostasis is largely preserved when noise is added to the input stimulus.

Although stochasticity in the input increases the seat slope, thereby somewhat degrading the homeostasis in $\langle y \rangle$, it also extends the range of $J$ values for which oscillations are produced, and in this way extends the interval of dynamic homeostasis. This occurs because, with stochastic input $J$, even when the mean $J$ is outside of the oscillation interval for the deterministic model, $J$ values from within the oscillation interval are still sampled. As a result, while the deterministic model would be at a stable equilibrium state, the model with stochastic input, but the same mean value of $J$, exhibits irregular oscillations and thereby extends the homeostatic range.

To quantify the extent at which the seat of the homeostasis curve was extended by the stochastic input, we computed the difference, for each $J$, between a stochastic chair curve in Fig. \ref{fgr7: Stochastic FHN model}A and the deterministic curve. Specifically, we looked for values of $J$ that pulled the stochastic chair curve away from the equilibrium branch of the deterministic chair curve. Once that difference was sufficiently large ($\Delta=0.03$), we used that $J$ value as an end point for the oscillation interval. This was done for each end point, and the distance between these two end points is plotted in Fig. \ref{fgr7: Stochastic FHN model}C as the ``oscillation interval length''. For each value of the feedback strength $\alpha$ the oscillation interval length is longer for the model with stochastic input than for the deterministic model. Indeed, the greater the noise strength the greater the extension of the oscillation interval.

\section{Bursting with the Chay-Keizer Model}
Teresa Chay and Joel Keizer published the first biophysical model for bursting electrical activity and Ca$^{2+}$ oscillations in pancreatic $\beta$-cells in the early 1980s \cite{chay1983minimal}, following the experimental work of Atwater et al. \cite{atwater1980nature}. This model, along with models for electrical bursting in neurons, was used by John Rinzel to carry out pioneering work on the analysis of bursting oscillations by decomposing the system into fast and slow subsystems \cite{rinzel1987formal,rinzel1986different}. This fast-slow analysis technique would be used in subsequent studies of bursting oscillators \cite{bertram1995topological,rinzel1987dissection,bertram2017multi,best2005dynamic}. In this section, we use the Chay-Keizer model and fast-slow analysis to demonstrate dynamic homeostasis in the slow variable that drives bursting oscillations. 

\subsection{Model formulation}
In the original form \cite{chay1983minimal}, the Chay-Keizer model is five-dimensional with differential equations for the membrane potential ($V$), two activation variables, one inactivation variable, and the free intracellular Ca$^{2+}$ concentration ($c$). However the important dynamics are preserved with a 3-dimensional reduction consisting of dynamics of $V$, the fraction of open K$^+$ channels ($w$), and $c$. The reduced model \cite{bertram2023deconstructing} is
\begin{align} 
    \frac{dV}{dt} &= -(I_{Ca}+I_K+I_{K(Ca)}+I_{K(ATP)}-I_{ap})/C_m, \label{eq4: Chay-Keizer voltage} \\
    \frac{dw}{dt} &= \frac{w_{\infty}(V) - w}{\tau_w}, \label{eq5: Chay-Keizer gating} \\
    \frac{dc}{dt} &= -f(\beta I_{Ca} + k_c c). \label{eq6: Chay-Keizer calcium}
\end{align}
The dynamics of $V$ (Eq. (\ref{eq4: Chay-Keizer voltage})) are governed by the different ionic currents and by the membrane capacitance $C_m$. The parameter $I_{ap}$ is the applied external current or a background current. The first term ($I_{Ca}$) is the inward Ca$^{2+}$ current that depolarizes the membrane and produces the upstroke of an action potential or spike. Its activation is fast and it is assumed to change instantaneously with $V$. The second term ($I_K$) is the outward delayed rectifier K$^+$ current, responsible for the downstroke of the action potential. Its activation occurs on a slower time scale ($\tau_w=16$ ms), and its dynamics are characterized by Eq. (\ref{eq5: Chay-Keizer gating}). The currents are given by:
\begin{align*} 
    I_{Ca} &= g_{Ca}m_{\infty}(V)(V-V_{Ca}), \\
    I_K &= g_K w(V-V_K),
\end{align*}
where $g_{Ca}$ and $g_K$ are maximum conductance parameters, $m_{\infty}$ is the equilibrium value for the fraction of open Ca$^{2+}$ channels, and $w$ is the fraction of open delayed rectifying K$^+$ channels. The constants $V_{Ca}$ and $V_K$ are the calcium and potassium reversal potentials, respectively. The equilibrium activation functions are $V$-dependent Boltzmann functions:
\begin{align*} 
    m_{\infty}(V) &= \frac{1}{1+\exp{(\frac{v_m - V}{s_m})}}, \\
    w_{\infty}(V) &= \frac{1}{1+\exp{(\frac{v_w - V}{s_w})}},
\end{align*}
where $v_m$ and $v_w$ are parameters that set the half-maximum response, and $s_m$ and $s_w$ set the slope of the activation curve.

\begin{table}[t!]
\centering
\begin{tabular}{c c c c} 
 \hline
 parameter & value & parameter & value \\ [0.5ex] 
 \hline
 $g_{Ca}$ & 1200 pS & $I_{ap}$ & 500 mV \\ 
 $g_{K(Ca)}$ & 300 pS & $p$ & 5 \\
 $g_K$ & 3000 pS & $K_{\Omega}$ & 0.3 $\mu$M \\ 
 $g_{K(ATP)}$ & 230 pS & $\tau_w$ & 16 ms  \\
 $C_m$ & 5300 fF & $\beta$ & $2.25\times10^{-6}$ $\mu$M fA$^{-1}$ ms$^{-1}$ \\
 $v_w$ & -16 mV & $s_w$ & 5 mV \\ 
 $v_m$ & -20 mV & $s_m$ & 12 mV \\
 $V_K$ & -75 mV & $f$ & 0.001 \\
 $V_{Ca}$ & 25 mV & $k_c$ & 0.07 ms$^{-1}$ \\
 \hline
\end{tabular}
\caption{Parameter values used in all figures for the Chay-Keizer model, except $k_c$ which is varied in several figures.  }
\label{table 1: Chay-Keizer parameters}
\end{table}

\begin{figure}[b!]
\centering
\includegraphics[width=12 cm]{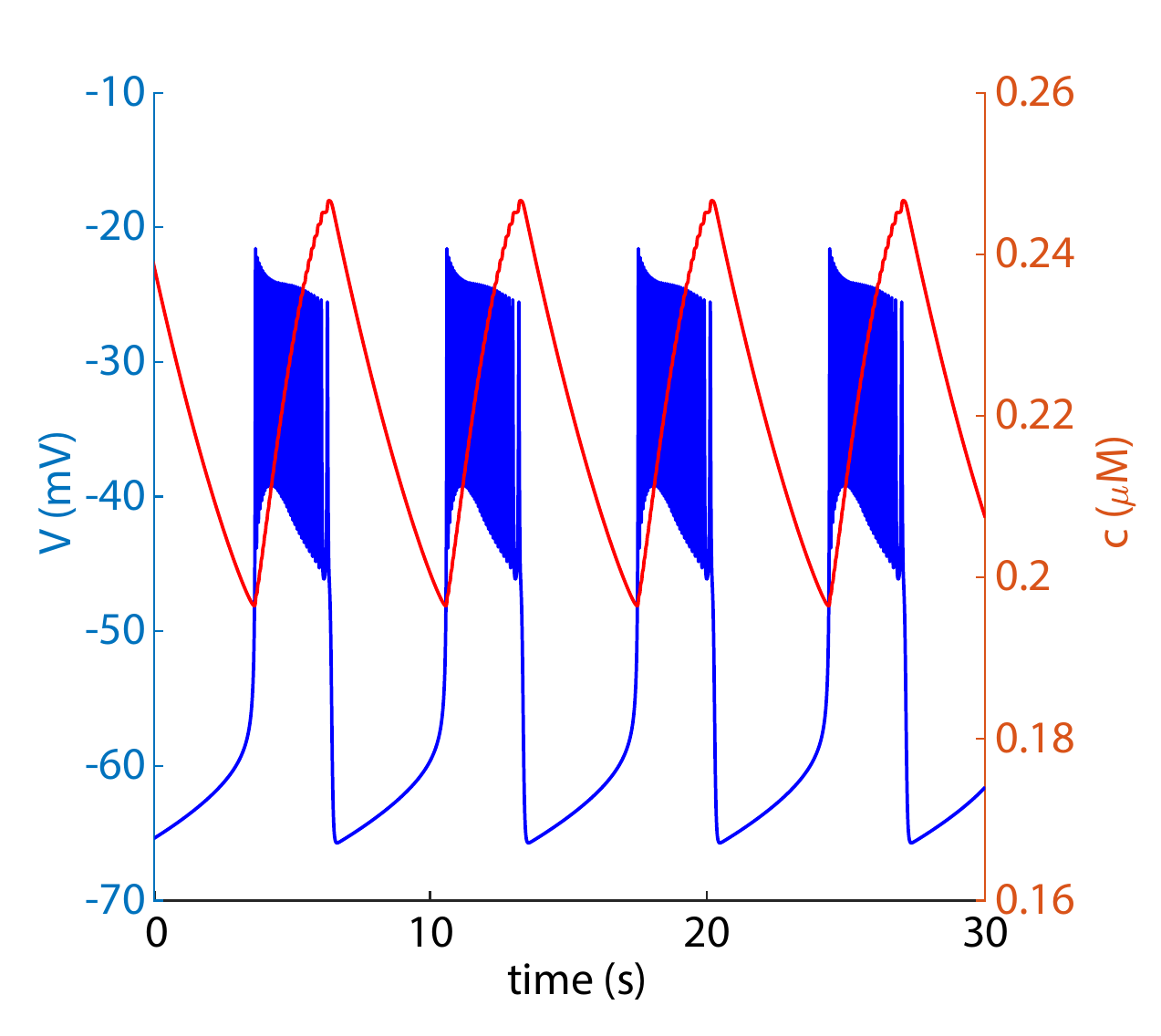}
\caption{Bursting oscillations produced by the reduced Chay-Keizer model. Bursting in membrane potential $V$ (blue curve), consisting of active phases where the cell is spiking and silent phases where the cell is hyperpolarized. The intracellular Ca$^{2+}$ concentration $c$ (orange curve) has a slow sawtooth oscillation pattern, rising during the active phase and falling during the silent phase. }
\label{fgr8: Chay-Keizer bursting oscillator}
\end{figure}

The third term ($I_{K(Ca)}$) in Eq. (\ref{eq4: Chay-Keizer voltage}) describes the ionic current through Ca$^{2+}$-activated K$^+$ channels in which the gating is determined by intracellular Ca$^{2+}$. The fraction of activated K(Ca) channels is described by a Hill function that increases with $c$, and the ionic current is given by:
\begin{equation*} 
    I_{K(Ca)} = g_{K(Ca)}(\frac{c^p}{K_{\Omega}^p + c^p})(V-V_K),
\end{equation*}
where $g_{K(Ca)}$ is the maximum conductance parameter, $p$ is the Hill coefficient, and the parameter $K_{\Omega}$ sets the intracellular Ca$^{2+}$ concentration for half-maximum activation of the channels.

The fourth term ($I_{K(ATP)}$) in Eq. (\ref{eq4: Chay-Keizer voltage}) is the ATP-sensitive K$^+$ current. Assuming constant conductance, the current satisfies
\begin{equation*}
    I_{K(ATP)} = g_{K(ATP)}(V-V_K),
\end{equation*}
where $g_{K(ATP)}$ is the conductance parameter.

The dynamics of $c$ (Eq. (\ref{eq6: Chay-Keizer calcium})) are characterized by the net flux of Ca$^{2+}$ across the membrane multiplied by $f$, the fraction of intracellular 
Ca$^{2+}$ that is free or unbound by buffers. The first term reflects Ca$^{2+}$ influx through Ca$^{2+}$ ion channels, so it is proportional to $I_{Ca}$. The second term reflects Ca$^{2+}$ removal through Ca$^{2+}$ pumps in the cell's membrane, and is assumed to be directly proportional to $c$. The parameter $\beta$ converts current to flux and the parameter $k_c$ is the Ca$^{2+}$ pump rate.

\subsection{Bursting oscillations and fast-slow analysis of the reduced Chay-Keizer model}
The reduced Chay-Keizer model produces electrical bursting with the parameter values listed in Table \ref{table 1: Chay-Keizer parameters}. Each burst is composed of an active phase of spiking and a silent phase where the membrane is hyperpolarized, shown in blue in Fig. \ref{fgr8: Chay-Keizer bursting oscillator}. At the beginning of an active phase, the Ca$^{2+}$ concentration $c$ builds up due to influx from voltage-dependent Ca$^{2+}$ channels. The rise in Ca$^{2+}$ activates K(Ca) current which terminates the spiking and ends a burst active phase.  The cell then enters a silent phase where there is a net removal of Ca$^{2+}$ by Ca$^{2+}$ pumps in the membrane. The Ca$^{2+}$ concentration and the K(Ca) conductance are slowly reduced. Once the K(Ca) conductance is sufficiently small, the membrane potential reaches spike threshold again, terminating the silent phase and starting a new active phase. 

\begin{figure}[h!]
\centering
\includegraphics[width=\columnwidth]{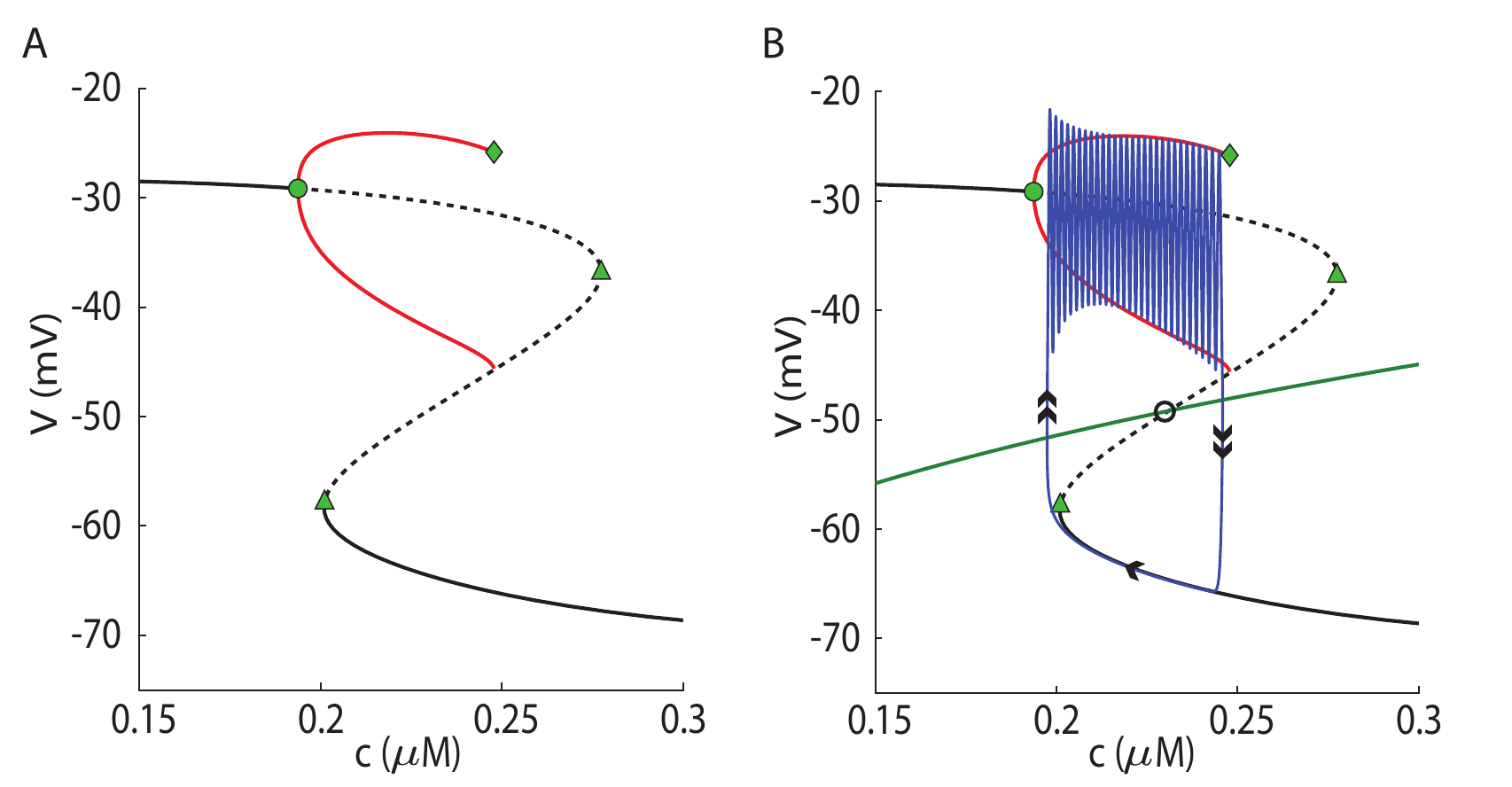}
\caption{Fast-slow analysis of the reduced Chay-Keizer model. (A) Bifurcation diagram of the fast subsystem. Stable (solid) and unstable (dashed) stationary equilibria are in black. The minimum and maximum $V$ of the stable periodic solutions are in red. There exists a supercritical Hopf bifurcation (green circle), two saddle-node bifurcations (green triangles), and a homoclinic bifurcation (green diamond). (B) The $c$-nullcline (green), full system equilibrium (black open circle), and burst trajectory (blue) are superimposed with the bifurcation diagram from panel A, now treated as a generalized $V$-nullcline. } 
\label{fgr9: Chay-Keizer fast/slow analysis}
\end{figure}

The fast-slow analysis of this bursting pattern relies on the slow dynamics of $c$ relative to the dynamics of the other two variables, $V$ and $w$, which form the fast subsystem.
A full description of the asymptotic dynamics of the fast subsystem over a range of $c$ values is shown in the bifurcation diagram in Fig. \ref{fgr9: Chay-Keizer fast/slow analysis}A, with $c$ as the bifurcation parameter. This is the ``critical manifold'' of the fast-slow system. 
Importantly, the fast subsystem exhibits bistability for values of $c$ between the leftmost saddle-node bifurcation (left green triangle) and the homoclinic bifurcation (green diamond), where both the periodic spiking solutions and low-$V$ equilibrium are stable. 

In Fig.~\ref{fgr9: Chay-Keizer fast/slow analysis}B we superimpose the $c$-nullcline (green curve) and the critical manifold in the $(c,V)$ plane. The intersection of the critical manifold and the $c$-nullcline (black open circle) is the equilibrium point of the full 3-dimensional model. Treating the critical manifold as a generalized $V$-nullcline, we use this and the $c$-nullcline to determine the flow of a trajectory in the ($c$,$V$)-plane, taking into account movement in the $c$ direction is much slower than movement in the $V$ direction due to the timescale separation. A burst trajectory, superimposed in blue, can be understood in terms of flow similar to the relaxation oscillation in Fig. \ref{fgr2: FHN Excitability and Oscillations}. During the silent phase, the trajectory slowly tracks along the bottom of the critical manifold until it reaches the left saddle-node or fold. After crossing this, the trajectory jumps up to the stable periodic (spiking) branch. Because in doing so the trajectory crosses the $c$-nullcline, it now flows to the right. The trajectory continues the spiking pattern while slowly moving rightward until it reaches the homoclinic bifurcation point. From here it jumps down to the stable stationary branch and switches direction to the left to begin a new silent phase.

\subsection{Dynamic homeostasis in the time-averaged calcium concentration} 
\begin{figure}[h!]
\centering
\includegraphics[width=\columnwidth]{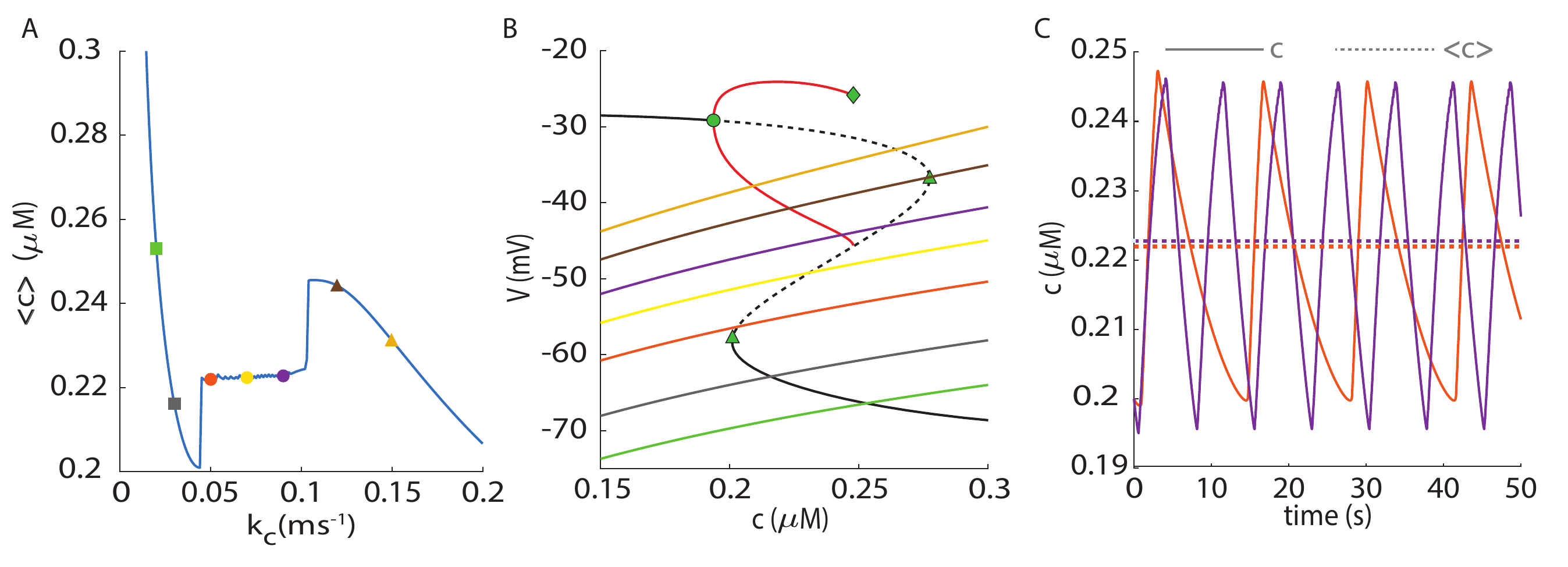}
\caption{Dynamic homeostasis in the time-average of the Ca$^{2+}$ concentration in the reduced Chay-Keizer model. (A) A homeostatic chair curve is evident over a range of Ca$^{2+}$ pump rate values. (B) Critical manifold and $c$ nullcline for several vales of $k_c$, color coded to match points in panel A. (C) Time series (solid) and time-averages (dashed lines) with $k_c=0.05$ ms$^{-1}$ (orange) and $k_c=0.09$ ms$^{-1}$ (purple) show  semi-invariance in the mean Ca$^{2+}$ concentration at equilibrium.}
\label{fgr10: Chay-Keizer slow variable invariance}
\end{figure}

We simulated the reduced Chay-Keizer model using Eqs. (\ref{eq4: Chay-Keizer voltage}-\ref{eq6: Chay-Keizer calcium}) with different values for the Ca$^{2+}$ pump rate ($k_c$). This is the parameter that is set by the glucose concentration; the pump rate is higher at a higher glucose concentration. We then calculated the asymptotic time-average of the slow variable, $\langle c \rangle$, with other parameters taken from Table \ref{table 1: Chay-Keizer parameters} (see Fig. \ref{fgr10: Chay-Keizer slow variable invariance}A). For very small or large values of $k_c$, increasing (decreasing) the pump rate results in a decrease (increase) in $\langle c \rangle$. Biologically, this makes sense because increasing (decreasing) the pump rate would increase (decrease) how much Ca$^{2+}$ is pumped out of a cell. However for intermediate values of $k_c$, changes in this pump rate result in very little change in the time-average of $c$. Thus, there is dynamic homeostasis and a chair-like curve for a bursting oscillator similar to the one in the previous section for a relaxation oscillator.

The seat of the chair curve in Fig. \ref{fgr10: Chay-Keizer slow variable invariance}A exists in the same parameter range of $k_c$ in which bursting oscillations occur. Fig. \ref{fgr10: Chay-Keizer slow variable invariance}B shows the critical manifold,
which does not change with $k_c$, and the $c$-nullcline, which does depend on $k_c$. These $c$-nullclines are color coded to correspond to the $k_c$ values at the points in Fig. \ref{fgr10: Chay-Keizer slow variable invariance}A. As $k_c$ is increased, the $c$-nullcline shifts upward. 

The chair curve can be partitioned into three sections that correspond to the three different asymptotic states of the model. For low values of $k_c$ (squares), the $c$-nullcline intersects the critical manifold on a stable stationary branch of the fast subsystem (Fig. \ref{fgr10: Chay-Keizer slow variable invariance}B). The full-system equilibrium at the intersection of these two curves is locally stable. Therefore, the membrane potential is characterized as 
fully quiescent at a low value of $V$. In this silent state, increasing the Ca$^{2+}$ pump rate lowers the mean Ca$^{2+}$ concentration, as can be seen in both panels A and B.

For intermediate values of $k_c$ (circles), the $c$-nullclines intersect the critical manifold at an unstable equilibrium point of the fast subsystem, so the full-system equilibrium is unstable. Also, in this parameter range the full system produces bursting oscillations, with dynamics similar to Fig. \ref{fgr9: Chay-Keizer fast/slow analysis}B. 
The time series for $c$ with $k_c$ set to 0.05 ms$^{-1}$ (orange) and 0.09 ms$^{-1}$ (purple) are shown in Fig. \ref{fgr10: Chay-Keizer slow variable invariance}C. The asymptotic time-averages of $c$ are superimposed as dashed horizontal lines. It is clear that the time-averages of the Ca$^{2+}$ concentration are nearly identical between the two time series, despite the oscillations having different duty cycles and periods. This semi-invariance occurs because the trajectory must cycle between the left saddle-node and homoclinic bifurcations of the critical manifold during each burst cycle, regardless of the oscillation period, and these bifurcations (and the entire critical manifold) are unaffected by $k_c$.


For higher values of $k_c$ (triangles), the $c$-nullcline intersects deep within the periodic branch of the critical manifold with the full system equilibrium remaining unstable. However, the trajectory no longer switches between the bottom stationary branch and periodic branch of the generalized $V$-nullcline to produce bursting. Instead, it  gets trapped within the periodic branch and the membrane potential continuously spikes \cite{terman1992transition}. Because the trajectory no longer cycles between the left saddle-node and homoclinic bifurcation, the invariance is lost, and the mean Ca$^{2+}$ concentration declines with larger values of the pump rate $k_c$. 

\subsubsection{The time-average of a fast variable is not homeostatic in a bursting oscillator}
We showed earlier that the time-average of the fast variable $x$ in the FHN model for relaxation oscillations is not homeostatic. 
To determine whether this is the case for a fast variable in the reduced Chay-Keizer model for bursting, we computed the time-average of the membrane potential, since this would be an observable variable in a $\beta$-cell (the other fast variable, $w$, cannot be directly measured). 
\begin{figure}[h!]
\centering
\includegraphics[width=\columnwidth]{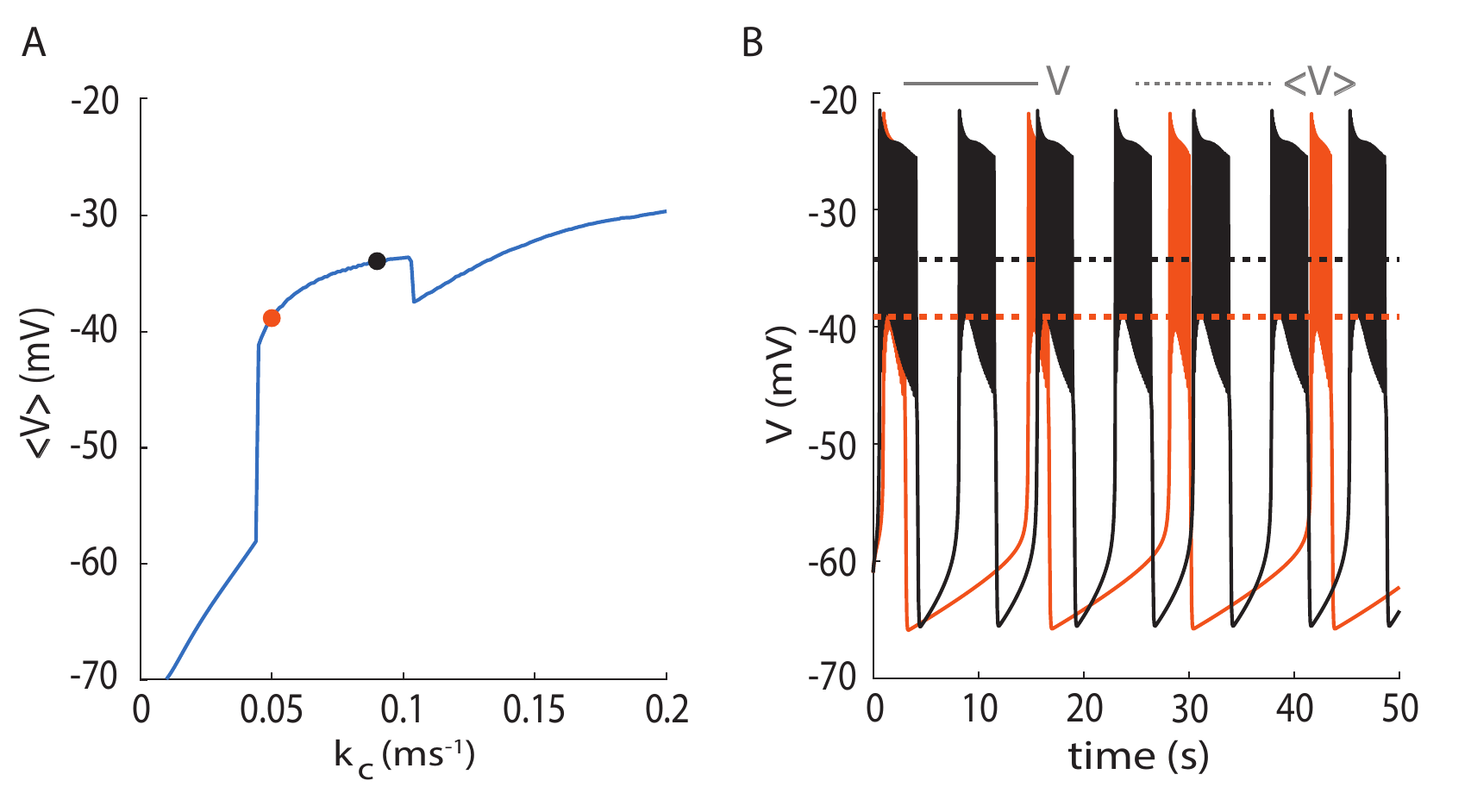}
\caption{Lack of dynamic homeostasis in the fast subsystem of the reduced Chay-Keizer model. (A) Time-average of $V$ plotted as a function of $k_c$ shows no dynamic homeostasis. (B) The time series (solid curves) for $k_c=0.05$ ms$^{-1}$ (orange) and $k_c=0.09$ ms$^{-1}$ (black) show that the burst duty cycle is much higher when $k_c$ is larger, raising the time-average of $V$ far above that for the lower $k_c$ value (dashed lines). } 
\label{fgr11: Chay-Keizer loss of invariance in fast variable}
\end{figure}
The time-average membrane potential, $\langle V \rangle$, as a function of $k_c$ is plotted in Fig. \ref{fgr11: Chay-Keizer loss of invariance in fast variable}A. Like before, the graph can be partitioned into three different segments corresponding to the three different dynamic states of the system: quiescence (low $k_c$), bursting (middle $k_c$), and continuous spiking (high $k_c$). This shows that the mean membrane potential is not invariant with respect to $k_c$, as with the fast variable of the relaxation oscillator.  

The time series of membrane potential when $k_c=0.05$ ms$^{-1}$ (orange) and $k_c=0.09$ ms$^{-1}$ (black) are shown in Fig. \ref{fgr11: Chay-Keizer loss of invariance in fast variable}B. The asymptotic time-averages of the membrane potential are superimposed as dashed lines. The duty cycle is lower for the lower $k_c$ value, driving the time-average down because the cell spends more time in the silent phase compared to the active phase. The duty cycle is higher with the higher $k_c$ value, producing a higher value of $\langle V \rangle$. The large difference in the duty cycle is responsible for the lack of dynamic homeostasis in the bursting fast subsystem. 

\subsection{Stochasticity in the Chay-Keizer model extends the range of homeodynamics}
Analogous to adding noise to the input $J$ in the FHN model, we now add noise to the glucose-dependent parameter $k_c$ in the reduced Chay-Keizer model:
\begin{figure}[b!]
\centering
\includegraphics[width=\columnwidth]{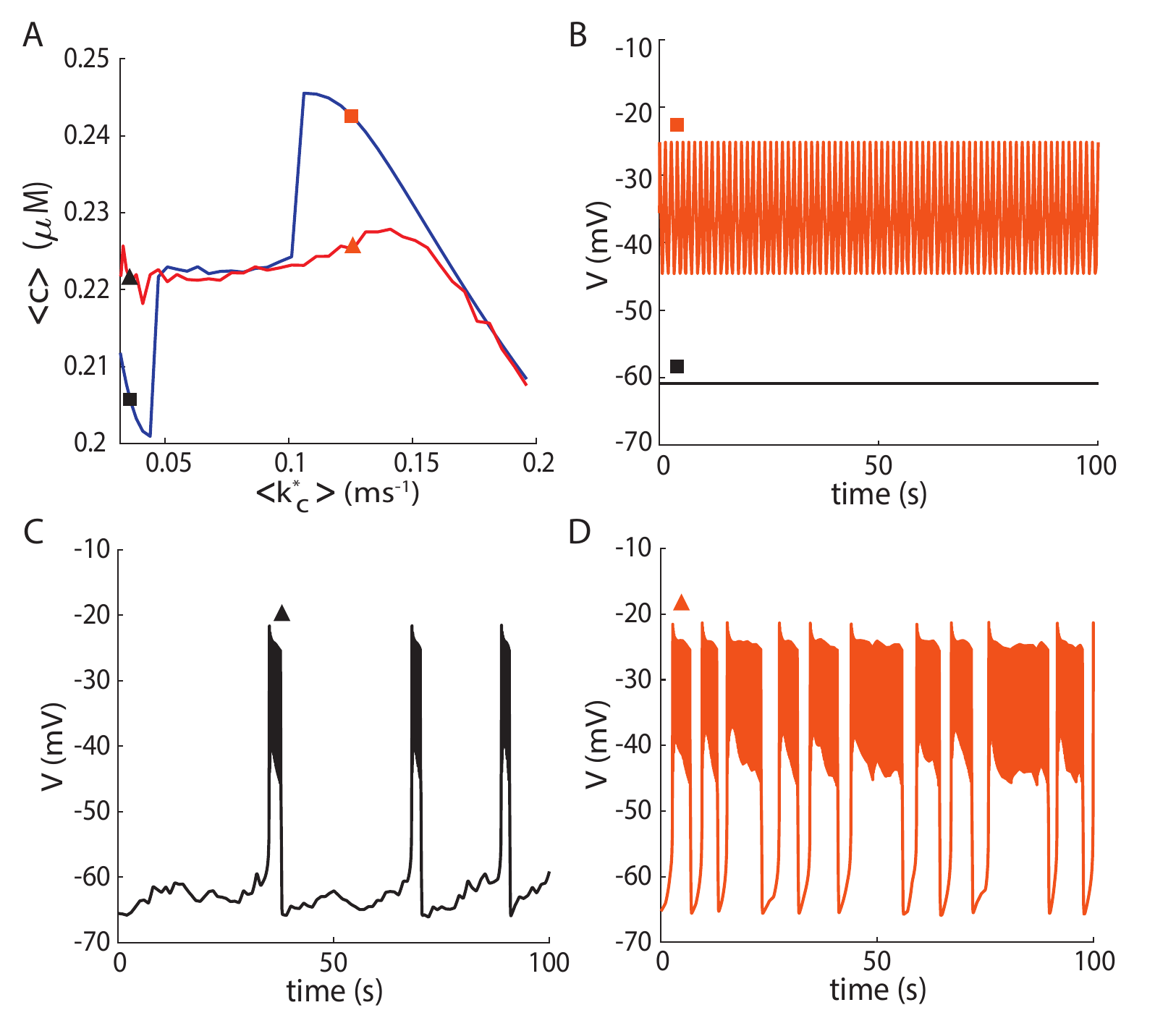}
\caption{Dynamic homeostasis is preserved and extended when stochasicity is added to the Ca$^{2+}$ pump rate of the reduced Chay-Keizer model. (A) Homeostatic chair curves for the deterministic (blue, $\sigma=0$) and stochastic (red, $\sigma=0.04$) versions of the model. (B) Membrane potential time series of the deterministic model when $k_c \approx 0.036$ (black square) and $k_c \approx 0.126$ (orange square) 
(C) Membrane potential time series of the stochastic model when $\langle k^*_c \rangle \approx 0.036$ (black triangle). (D) Membrane potential time series of stochastic model when $\langle k^*_c \rangle \approx 0.126$ (orange triangle). } 
\label{fgr12: Stochastic Chay-Keizer model}
\end{figure}
\begin{align} 
    \frac{dV}{dt} &= -(I_{Ca}+I_K+I_{K(Ca)}+I_{K(ATP)}-I_{ap})/C_m \label{eq7: Stochastic Chay-Keizer voltage} \\
    \frac{dw}{dt} &= \frac{w_{\infty}(V) - w}{\tau_w} \label{eq8: Stochastic Chay-Keizer gating} \\
    \frac{dc}{dt} &= -f(\beta I_{Ca} + k_c^* c). \label{eq9: Stochastic Chay-Keizer calcium}
\end{align}
where all variables and parameters are defined as in (Eqs. (\ref{eq4: Chay-Keizer voltage}-\ref{eq6: Chay-Keizer calcium})), except the new calcium pump rate $k_c^*$, which now follows a folded normal distribution:
\begin{equation} \label{eq10: Folded normal distribution}
    k_c^* \sim FN(k_c,\sigma).
\end{equation}
This distribution, a variant of the normal distribution, was used to ensure that the pump rate remains non-negative. Details of the folded normal distribution are given in the Appendix \ref{Appendix: Folded Normal distribution}.

We numerically solved the system  (\ref{eq7: Stochastic Chay-Keizer voltage}-\ref{eq9: Stochastic Chay-Keizer calcium}) using the forward Euler method, randomly assigning a new calcium pump rate every 1 second according to (\ref{eq10: Folded normal distribution}). We then calculated the time-average of the Ca$^{2+}$ concentration, $\langle c \rangle$, after sufficiently long time for many different values of the folded normal distribution mean, $k_c$. This is plotted in Fig. \ref{fgr12: Stochastic Chay-Keizer model}A as a function of the expected value of the pump rate, $\langle k^*_c \rangle$ (see Appendix \ref{Appendix: Folded Normal distribution} for the expected value formula for the folded Normal distribution). The blue curve corresponds to the deterministic simulation (when $\sigma=0$), while the red curve is the stochastic simulation (when $\sigma=0.04$).

For low values of $\langle k^*_c \rangle$, 
$\langle c \rangle$ for the stochastic model is relatively invariant and has a higher value than in the deterministic case. The reason for this is evident from the $V$ time courses of the stochastic and deterministic model. For a low value of 
$k_c$ in the deterministic model (denoted by a black square in Fig. \ref{fgr12: Stochastic Chay-Keizer model}A), the membrane potential is quiescent at a low value (black curve in Fig. \ref{fgr12: Stochastic Chay-Keizer model}B). However, at the same mean value of $k_c^*$ in the stochastic model (black triangle), $k_c^*$ is sampled at some values in which bursting is produced, resulting in an irregular burst pattern (Fig. \ref{fgr12: Stochastic Chay-Keizer model}C). This results in higher Ca$^{2+}$ time-averages.

For intermediate values of $\langle k^*_c \rangle$, the stochastic and deterministic Ca$^{2+}$ time-averages are both semi-invariant, forming the seat of the chair curve. Within this range, both the deterministic and stochastic time series exhibit regular bursting (not shown). 

For high value of $\langle k^*_c \rangle$, the seat of the chair extends considerably farther in the stochastic model than it did in the deterministic model. The reason for this is again evident from an examination of the time series. At a large $k_c$ value, the deterministic model (orange square)  produces continuous spiking (orange curve in Fig. \ref{fgr12: Stochastic Chay-Keizer model}B), a behavior in which dynamic homeostasis does not occur. For the stochastic model with the same mean value of $k_c$ (orange triangle), the membrane potential produces periods of regular bursting with moments of extended silent and active phases (see Fig. \ref{fgr12: Stochastic Chay-Keizer model}D). Because the burst pattern is dominant, dynamic homeostasis occurs at and near this value of $\langle k^*_c \rangle$. Also, due to the silent phases between active spiking phases, the mean Ca$^{2+}$ concentration is lower than in the deterministic model. 

\section{Dynamic Homeostasis Occurs in Bursting Oscillations with Multiple Slow Variables}
\begin{figure}[b!]
\centering
\includegraphics[width=\columnwidth]{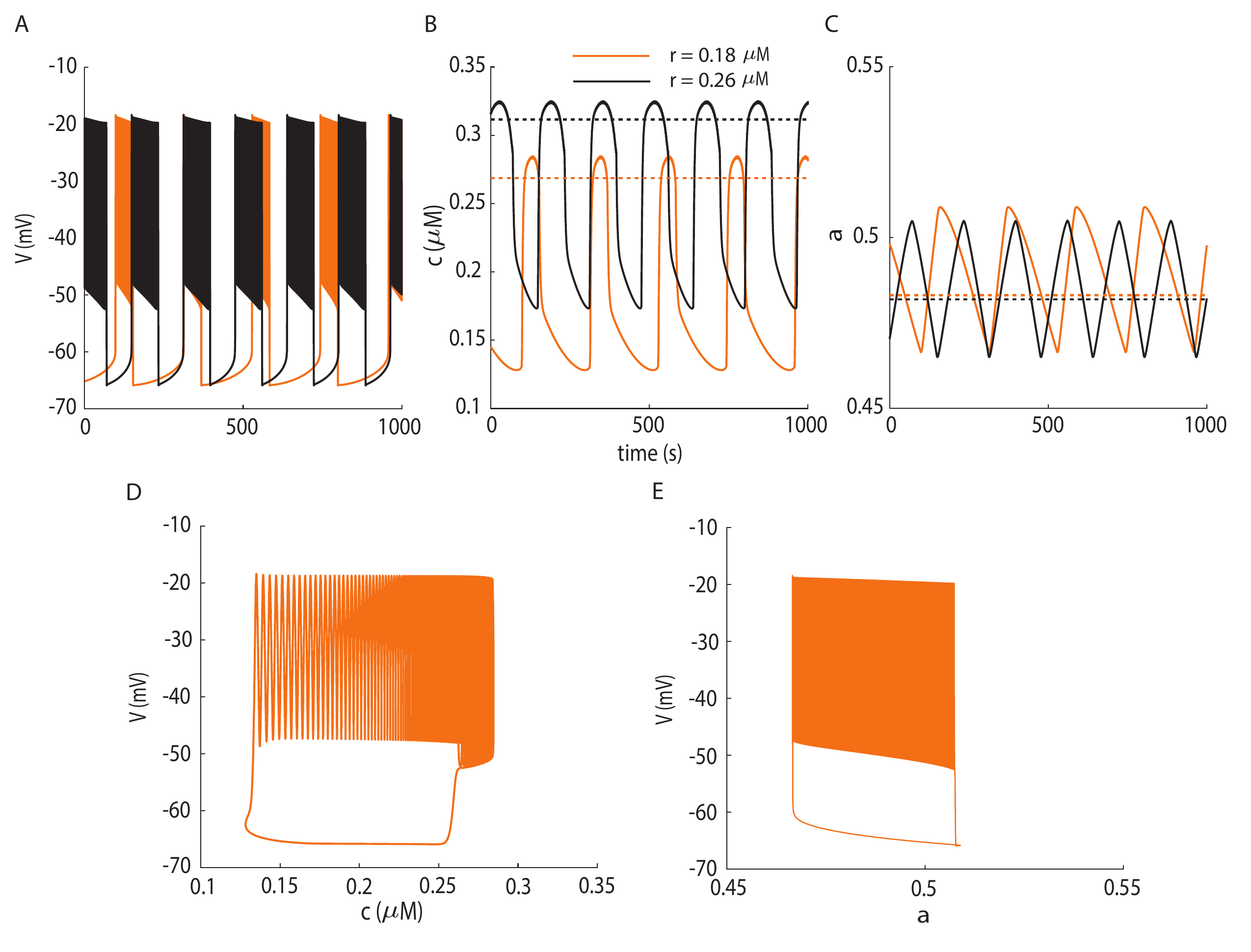}
\caption{Dynamic homeostasis in the mean $\frac{[ADP]}{[ATP]}$ ratio $a$ in the PBM model when $g_{K(Ca)}=25$ pS and bursting is slow. Time series for (A) membrane potential, (B) cytosolic Ca$^{2+}$ concentration, and (C) $\frac{[ADP]}{[ATP]}$ ratio when $r=0.18$ $\mu$M (orange) and $r=0.26$ $\mu$M (black). The time-average of $a$ exhibits dynamic homeostasis, but the time-average of $c$ does not. (D) The burst trajectory in ($c,V$)-space moves to the right during the first part of the burst active phase, and moves back to the left during the second part. Only the case with $r=0.18$ $\mu$M is shown here and in the next panel. (E) The burst trajectory in $(a,V)$-space, with $a$ moving monotonically to the right during the burst active phase, and to the left during the burst silent phase. Other parameters are given in Table \ref{table 2: Phantom burster parameters}.} 
\label{fgr13: Phantom Burster gkca Low}
\end{figure}
We have shown that dynamic homeostasis can be achieved in the FitzHugh-Nagumo model where a single slow variable drives relaxation oscillations, and in the reduced Chay-Keizer model in which a single slow variable (the Ca$^{2+}$ concentration) drives the bursting oscillations. What happens when there are multiple slow variables? Does dynamic homeostasis still occur, and if so, which of the slow variables maintains semi-invariance in the mean as an input parameter is varied? To address these questions, we studied the phantom bursting model (PBM) for pancreatic $\beta$-cells \cite{bertram2000phantom,bertram2004calcium}.

The PBM is an extension of the reduced Chay-Keizer model with additional dynamics for the Ca$^{2+}$ concentration in the endoplasmic reticulum (ER, denoted as $c_{er}$) and the ratio of adenosine diphosphate (ADP) to adenosine triphosphate (ATP), $a=\frac{[ADP]}{[ATP]}$. The differential equations for this ``Phantom Bursting Model'' (PBM) are: 

\begin{equation}\label{eq11: Phantom bursting model}
    \begin{split}
        \frac{dV}{dt} &= -(I_{Ca}+I_K+I_{K(Ca)}+I_{K(ATP)})/C_m, \\
        \frac{dw}{dt} &= \frac{w_{\infty}(V)-w}{\tau_w}, \\
        \frac{dc}{dt} &= f_{cyt}(J_{mem}+J_{er}), \\
        \frac{c_{er}}{dt} &= -f_{er}(V_{cyt}/V_{er})J_{er}, \\
        \frac{da}{dt} &= \frac{a_{\infty}(c)-a}{\tau_a}.
    \end{split}
\end{equation}
The details for the model variables and parameters are discussed in Appendix \ref{Appendix: Phantom Bursting}. System (\ref{eq11: Phantom bursting model}) includes 3 slow variables ($c$, $c_{er}$, and $a$), and each can independently drive bursting oscillations in the membrane potential depending on the parameters. Of these three, the cytosolic Ca$^{2+}$ concentration $c$ is the fastest, with the ER Ca$^{2+}$ concentration $c_{er}$ varying on a slower time scale. The nucleotide ratio, $a$, varies on the slowest time scale and is itself a function of $c$ rather than voltage $V$. These variables are discussed in detail in \cite{bertram2004calcium}, and a detailed time scale analysis of $V$, $c$, and $c_{er}$ was performed in \cite{mckenna18}. The time scale for $a$ is simply the constant $\tau_a=5$ min. It was shown in \cite{bertram2004calcium} that when the conductance of the K(Ca) current ($g_{K(Ca)}$) is large any bursting that occurs is fast, with a period of approximately 10 s. This fast bursting is driven entirely or almost entirely by the dynamics of $c$. If the K(Ca) conductance is small, then the bursting that occurs is slow, with a period of approximately 5 min, and is driven primarily by the dynamics of $a$.  

\begin{figure}[b!]
\centering
\includegraphics[width=\columnwidth]{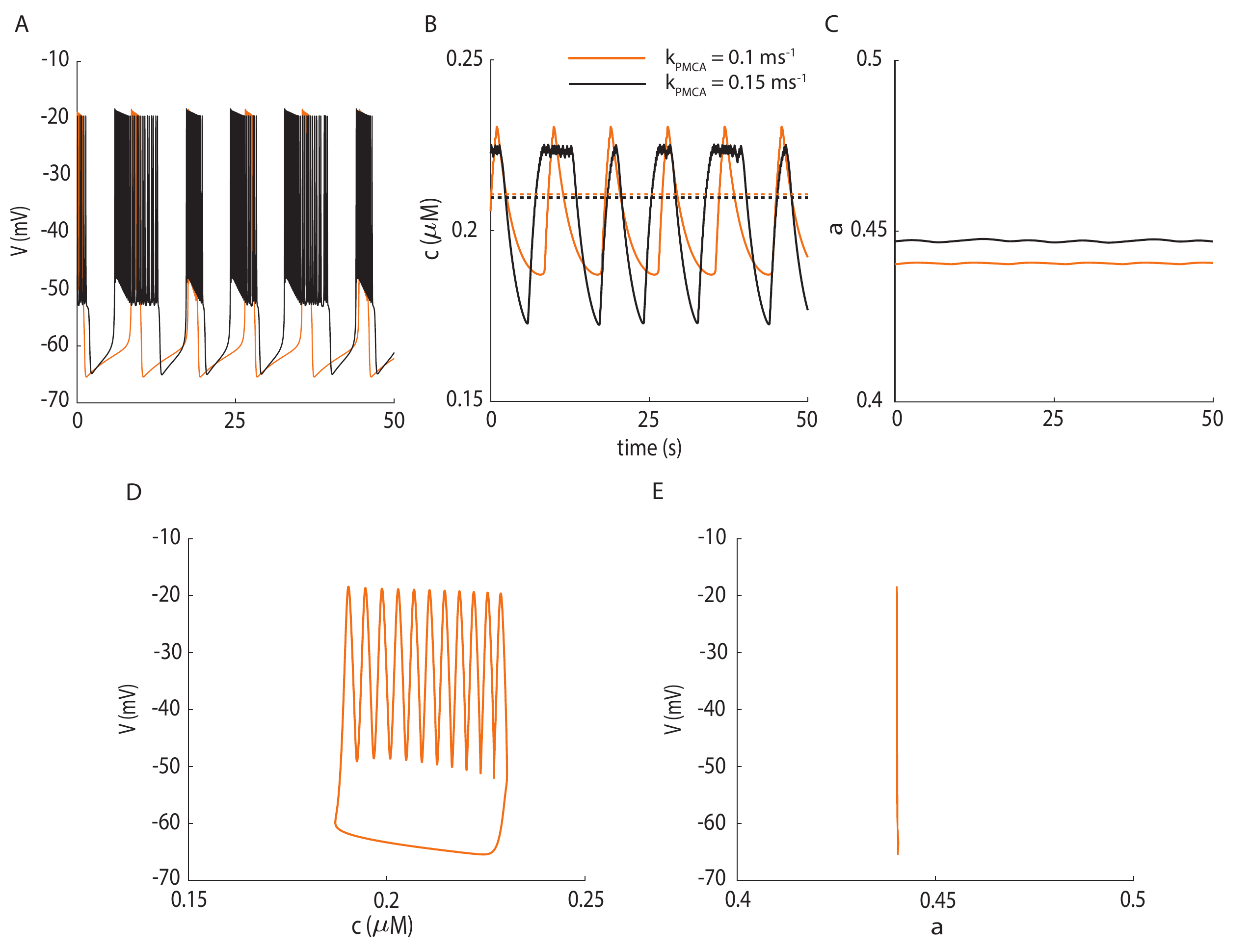}
\caption{Dynamic homeostasis in the mean cytosolic Ca$^{2+}$ concentration $c$ in the PBM model when $g_{K(Ca)}=600$ pS. Time series for (A) membrane potential, (B) cytosolic Ca$^{2+}$ concentration, and (C) $\frac{[ADP]}{[ATP]}$ ratio when $k_{PMCA}=0.1$ ms$^{-1}$ (orange) and $k_{PMCA}=0.15$ ms$^{-1}$ (black). (D) The burst trajectory projected into the ($c,V$)-plane moves to the right during the active phase and to the left during the silent phase, having a profile typical for a variable that drives bursting. Only the case with $k_{PMCA}=0.1$ ms$^{-1}$ is shown here and in the next panel. (E) The projection of the burst trajectory into the $(a,V)$-plane appears as a vertical line segment since $a$ is relatively constant during this fast bursting and plays no role in driving the bursting. Other parameters are given in Table \ref{table 2: Phantom burster parameters}.} 
\label{fgr14: Phantom Burster gkca High}
\end{figure}

Fig. \ref{fgr13: Phantom Burster gkca Low} shows the dynamics of $V$, $c$, and $a$ when the Ca$^{2+}$-activated $K^+$ maximum conductance is at a low value, $g_{K(Ca)}=25$ pS. To determine whether this system can display dynamic homeostasis, we varied the metabolism parameter $r$, which is a parameter in the $a$ differential equation (see Appendix C). Increasing $r$ translates the $a$-nullcline upward in the $(a,V)$ plane (not shown), playing a role similar to that played by $k_c$ with $c$ in the reduced Chay-Keizer model. When $r$ is increased from 0.18 $\mu M$ (orange curves) to 0.26 $\mu M$ (black curves), there is a decrease in the burst period and an increase in its duty cycle (Fig. \ref{fgr13: Phantom Burster gkca Low}A), and for this reason $r$ is considered the glucose-sensing parameter in the model. These changes have a big impact on $c$, which has a smaller amplitude and considerably larger time-average with the higher $r$ value (Fig. \ref{fgr13: Phantom Burster gkca Low}B). In contrast, the time-average of $a$ is almost invariant (Fig. \ref{fgr13: Phantom Burster gkca Low}C), indicating dynamic homeostasis in the time-average of this slow variable to changes in the glucose parameter $r$. 

As illustrated in Fig. \ref{fgr10: Chay-Keizer slow variable invariance}, the dynamic homeostasis occurs because the slow variable that drives bursting, i.e., starts and stops active phases of spiking, oscillates over the same range of values when a parameter, like $k_c$ in that case and $r$ in this case, is varied. The active phase starts at a fast-subsystem saddle-node bifurcation and ends when the trajectory moves past a fast-subsystem homoclinic bifurcation. All that changes is the location of the slow nullcline. We plotted the solution trajectory when $r=0.18$ $\mu M$ in the $(c,V)$-plane in Fig. \ref{fgr13: Phantom Burster gkca Low}D. During the silent and much of the active phase the burst trajectory looks much like it did in Fig. \ref{fgr9: Chay-Keizer fast/slow analysis}B. However, rather than crossing through the fast-subsystem homoclinic bifurcation, the trajectory stalls, and even moves back to the left before the active phase terminates. This indicates that the K(Ca) ionic current activated by $c$ is insufficient to terminate the active phase. It also indicates that $c$ is not driving the bursting. The variable that is driving the bursting is $a$. When the trajectory is viewed in the $(a,V)$-plane, it looks just as it did in Fig. \ref{fgr9: Chay-Keizer fast/slow analysis}B, but with a different slow variable. That is, the active phase starts at a fast-subsystem saddle-node bifurcation and ends at a fast-subsystem homoclinic bifurcation (the fast subsystem would be the variables other than $a$); $a$ is driving the bursting in this case. This projection of the trajectory looks the same for the two values of parameter $r$, and explains the dynamic homeostasis in the time-average of $a$.

Fig. \ref{fgr14: Phantom Burster gkca High} shows the dynamics of $V$, $c$, and $a$ when the Ca$^{2+}$-activated $K^+$ maximum conductance is large, $g_{K(Ca)}=600$ pS. Since bursting is now driven by $c$, we varied a parameter in the $c$ differential equation, the Ca$^{2+}$ pump rate $k_{PMCA}$ (called $k_c$ in the Chay-Keizer model). When $k_{PMCA}$ was raised from 0.1 ms$^{-1}$ (orange curves) to 0.15 ms$^{-1}$ (black curves), the burst duty cycle increased, similar to before in the Chay-Keizer model. Unlike in Fig. \ref{fgr13: Phantom Burster gkca Low}, dynamic homeostasis was present in the time-average of $c$ (see Fig. \ref{fgr14: Phantom Burster gkca High}B). The time-average of $a$ was not invariant, but increased with $k_{PMCA}$ (see Fig. \ref{fgr14: Phantom Burster gkca High}C). From the projection of the burst trajectory into the $(c,V)$-plane (Fig. \ref{fgr14: Phantom Burster gkca High}D), it is apparent that the bursting is driven by $c$, just as with the Chay-Keizer model in Fig. \ref{fgr9: Chay-Keizer fast/slow analysis}B. In contrast, the trajectory projected into the $(a,V)$-plane (Fig. \ref{fgr14: Phantom Burster gkca High}E) shows that $a$ is almost constant during this fast bursting and plays no role in driving the burst.  

The central conclusion from this analysis is that dynamic homeostasis is present in a bursting model with multiple slow variables, but only in the slow variable that drives the bursting oscillations.

\section{Discussion}
Homeostasis is a fundamental principle in biology, and much recent work has been done on understanding it from a dynamical systems viewpoint in signaling pathways \cite{reed2017analysis,nijhout2014escape,nijhout2019systems}. The focus there was on the semi-invariance of a system variable at static equilibrium in response to change in input. Other works \cite{lloyd2001homeodynamics,xiong2023physiological} have argued for the reemergence of a more general approach to homeostasis, the concept of dynamic homeostasis. Real biological systems are complex with many degrees of freedom allowing for more complicated outputs like bistability, oscillations, chaos, or possibly some combination. Everyone has heard how normal body temperature is 98.6$^\circ$F (37$^\circ$C) or that normal blood glucose level is about 85 mg/dL (4.7 mmol/L). However, these are only averages. Body temperature \cite{aschoff1971human}, glucose, and insulin \cite{sturis1991computer} all oscillate over time even with constant, controlled external stimuli present. 

In this article, we demonstrated that homeostasis can occur even in oscillatory systems. In particular, this dynamic homeostasis occurs in oscillatory systems with fast-slow dynamics such as relaxation oscillations and electrical bursting oscillations. We also demonstrated that the homeostasis is restricted to the slow variable that drives the oscillation; fast variables do not exhibit the behavior, nor do additional slow variables not responsible for driving oscillations. The dynamic homeostasis is robust to noise in the input to the slow subsystem, and in fact this type of noise can extend the range of input at which dynamic homeostasis occurs.



What advantage might there be to a cell to maintain a constant average level of a slow variable over a range of inputs? Of course, this depends on the slow variable. For the Chay-Keizer model, there is semi-invariance in the mean cytosolic Ca$^{2+}$. This could be desirable as a means to ensure that the Ca$^{2+}$ level does not go too high even as the variable's production is increased. This ion regulates many enzymes and transcription factors, and some ion channels. It is also a trigger for secretion, as we discuss below. Programmed cell death, apoptosis, is also triggered by sustained Ca$^{2+}$ levels in the cytoplasm and mitochondria of a cell. Therefore, the intracellular Ca$^{2+}$ level should be tightly regulated, and dynamic homeostasis as in the Chay-Keizer model is one way to achieve this. Likewise, creatine kinase is found in many tissues and acts as an ATP buffer, storing ATP for future use if needed. The dynamic homeostasis of the mean $\frac{[ADP]}{[ATP]}$ level is another means for keeping the nucleotide ratio constant even as the production of ATP from ADP is increased. In ecological systems, the birth/death dynamics of predator and prey can vary substantially, giving rise to population dynamics with fast and slow subsystems. This is true, for example, in resource-consumer interactions, where the life cycle of the consumer (e.g., insect herbivores) can be much shorter than that of the resource (the trees that they feed on). Such systems can exhibit relaxation oscillations \cite{hastings18}, and while the resource level may oscillate over time, its average value would be homeostatic over a large range of parameter values in its subsystem.

The two models for bursting used in this article were developed to describe bursting in $\beta$-cells from pancreatic islets of Langerhans. These cells secrete insulin in response to the blood glucose level; higher glucose levels, as would occur after a meal, result in higher insulin secretion. Therefore, the primary input to these cells is glucose and the primary output is insulin. The transduction pathway between input and output is very complicated \cite{rorsman18,satin15}, and has been the focus of mathematical modeling since 1983 with the publication of the Chay-Keizer model \cite{chay1983minimal}. In that original model, the input (glucose) acted directly on the Ca$^{2+}$ pump rate, $k_c$, since the product of glucose metabolism in the cell, ATP, powers those pumps. Although a very sensible notion, we have shown here that with that model and that source of input, the temporal average level of Ca$^{2+}$ is almost constant over the entire range of input values for which the cell is bursting (the state of the cell when the animal is in a fed state). Since insulin secretion is evoked by the binding of Ca$^{2+}$ to proteins in the exocytotic machinery, this would lead to a relatively steady level of insulin secretion at all stimulatory glucose levels if the relationship between Ca$^{2+}$ concentration and insulin were linear. The reality is much more complicated than that \cite{sedaghat02}, but an invariant average Ca$^{2+}$ level with changes in glucose is nonetheless a very undesirable feature, and is a byproduct of the key role that the Ca$^{2+}$ concentration plays in driving bursting in that model.

In the Phantom Bursting Model used in the latter part of this study, the role of burst driver during the slow bursting typically observed in islets is played by a different variable, $a$, the ratio of ADP to ATP in the cell. This ratio sets the activation level of ATP-sensitive K$^+$ channels, K(ATP) channels, which were not discovered until after the publication of the Chay-Keizer model \cite{ashcroft84,cook84b,henquin88b}. Since $a=\frac{[\rm{ADP}]}{[\rm{ATP}]}$ now drives the bursting, the cytosolic Ca$^{2+}$ concentration is released from homeostasis, so it can, and does, increase with higher glucose concentrations. This is much more appropriate for $\beta$-cell function. Thus, dynamic homeostasis of the nucleotide ratio may simply be a byproduct of a design that aims to have increased insulin secretion with higher levels of blood glucose. The semi-invariance in the nucleotide ratio in $\beta$-cell models in which bursting is driven by the nucleotide ratio was first described in \cite{marinelli2022oscillations}, and was supported by experimental data in the same publication. This fact was used to make the case that slow bursting in these cells was driven by oscillations in the K(ATP) channel conductance, rather than through K(Ca) conductance. A much earlier study recognized the semi-invariance of the cytosolic Ca$^{2+}$ concentration in a model where bursting was driven by Ca$^{2+}$ through actions on K(Ca) channels \cite{himmel87}. 

What this study adds to these prior studies is a demonstration of the persistence of this phenomenon across oscillators with a fast-slow structure, and the robustness of the phenomenon to noise. Also, the manner in which the dynamic homeostasis manifests in systems with multiple slow variables was described here for the first time. In particular, we demonstrated that it is only the slow variable responsible for driving bursting that exhibits dynamic homeostasis. This generalizes naturally to relaxation oscillators with multiple slow variables.


\appendix
\section{Linear stability analysis and Hopf bifurcation for the FHN model}\label{Appendix: FHN Linear stability and Hopf bifurcation}
In this appendix, we linearize the FHN model to find equations for Hopf bifurcations of the FHN model. The elements of the vector field for the FHN model are
\begin{equation*} 
    g(x,y) =\mu(x-\frac{x^3}{3}-y) \qquad h(x,y) = \frac{1}{\mu}(J+\alpha x-y),
\end{equation*}
where an equilibrium $(x^*,y^*)$ satisfies $g(x^*,y^*)=h(x^*,y^*)=0$, or
\begin{equation*}
    \frac{{x^*}^3}{3} + (\alpha-1)x^*+J=0 \qquad y^*=J+\alpha x^*.
\end{equation*}
The Jacobian for the linearized system is
\begin{equation*}
    A= \begin{pmatrix}
    \mu(1-{x^*}^2) & -\mu \\ \frac{\alpha}{\mu} & -\frac{1}{\mu}
    \end{pmatrix},
\end{equation*}
where the trace and determinant are $T = \mu(1-{x^*}^2)-\frac{1}{\mu}$ and $D = {x^*}^2+\alpha-1$.  A saddle-node bifurcation occurs when there is one zero eigenvalue and one non-zero eigenvalue, so $D=0$ and $T \neq 0$. This occurs at $\alpha=1$ and $x^*=0$. For $\alpha < 1$ there are three equilibria, while for $\alpha > 1$ there is a single equilibrium. We consider only the latter case in this study. 
A Hopf bifurcation 
occurs when the eigenvalues of the Jacobian are purely imaginary, so $T=0$ and $D\neq0$. These conditions give the locations and conditions for a Hopf bifurcation in terms of the parameters
\begin{equation} \label{eq12: Hopf conditions}
\begin{split}
x_{\pm}^*&=\pm\sqrt{1-\frac{1}{\mu^2}}, \\
J_{\pm}&=\pm(1-\alpha)(1-\frac{1}{\mu^2})^{\frac{1}{2}}\mp\frac{1}{3}(1-\frac{1}{\mu^2})^{\frac{3}{2}}, \\
\alpha&>\frac{1}{\mu^2}.
\end{split}
\end{equation}
In this manuscript, we assume $\mu\gg1$ to ensure a large timescale separation between the $x$ and $y$ species, and since we only consider the case $\alpha > 1$ the last condition is always satisfied. Given $\mu\gg1$ we can approximate the Hopf bifurcation curves in (\ref{eq12: Hopf conditions}) as $J_{\pm} \approx \pm(\frac{2}{3}-\alpha)$. The Hopf bifurcation curves are shown in a 2-parameter bifurcation diagram in Fig. \ref{fgr6: FHN invariance dependance on alpha}C. For $\frac{2}{3}-\alpha<J<\alpha-\frac{2}{3}$, the equilibrium is unstable and stable periodic solutions can exist. Otherwise, the equilibrium is locally stable and no periodic solutions exist.

\section{The Folded Normal Distribution}\label{Appendix: Folded Normal distribution}
The folded normal distribution was used in the stochastic Chay-Keizer model to ensure the Ca$^{2+}$ pump rate, $k_c$, remained non-negative during simulations. Given a normally distributed random variable $X$ with mean $k_c$ and standard deviation $\sigma$, the random variable $k_c^*=|X|$ has a folded normal distribution. The transformed expected value and standard deviation are as follows
\begin{equation} \label{eq13: Folded normal distribution}
\begin{split}
\langle k_c^* \rangle &= \sqrt{\frac{2}{\pi}} \sigma \exp(-\frac{k_c^2}{2\sigma^2}) + \mu(1-2\Phi(-\frac{k_c}{\sigma})),\\
\sigma_f &= \sqrt{k_c^2+\sigma^2 - \langle k_c^* \rangle^2},
\end{split}
\end{equation}
where $\Phi(x) = \frac{1}{2}(1-\text{erf}(\frac{x}{\sqrt{2}}))$ is the normal cumulative distribution function. The function $\text{erf}(z)=\frac{2}{\sqrt{\pi}} \int_{0}^{z} e^{-t^2}dt$ is the Gauss error function.

\section{The Phantom Bursting Model}\label{Appendix: Phantom Bursting} 
The phantom bursting model (PBM) is a pancreatic cell model first introduced and analyzed in \cite{bertram2000phantom}. The model includes the dynamics of 5 variables as follows:
\begin{equation*}
    \begin{split}
        \frac{dV}{dt} &= -(I_{Ca}+I_K+I_{K(Ca)}+I_{K(ATP)})/C_m, \\
        \frac{dw}{dt} &= \frac{w_{\infty}(V)-w}{\tau_w}, \\
        \frac{dc}{dt} &= f_{cyt}(J_{mem}+J_{er}), \\
        \frac{c_{er}}{dt} &= -f_{er}(V_{cyt}/V_{er})J_{er}, \\
        \frac{da}{dt} &= \frac{a_{\infty}(c)-a}{\tau_a}.
    \end{split}
\end{equation*}
where $V$ is the membrane potential, $w$ is the delayed rectifier activation, $c$ is the free cytosolic Ca$^{2+}$ concentration, $c_{er}$ is the free Ca$^{2+}$ concentration in the endoplasmic reticulum (ER), and $a$ is the nucleotide ratio ($\frac{[ADP]}{[ATP]}$). The parameters used in this model are given in Table \ref{table 2: Phantom burster parameters}.
\begin{table}[h!]
\centering
\begin{tabular}{c c c c} 
 \hline
 parameter & value & parameter & value \\ [0.5ex] 
 \hline
 $g_{Ca}$ & 1200 pS & $r$ & 0.225 $\mu$M \\ 
 $g_{K(Ca)}$ & 600 pS & $V_{cyt}/V_{er}$ & 10 \\
 $g_K$ & 3000 pS & $K_d$ & 0.4 $\mu$M \\ 
 $g_{K(ATP)}$ & 500 pS & $\tau_w$ & 18 ms  \\
 $C_m$ & 5300 fF & $\tau_a$ & 300,000 ms \\
 $v_w$ & -15 mV & $s_w$ & 5 mV \\ 
 $v_m$ & -20 mV & $s_m$ & 12 mV \\
 $p_{leak}$ & 0.0002 ms$^{-1}$ & $s_a$ & 0.1 $\mu$M\\
 $V_K$ & -75 mV & $f_{cyt}$ & 0.001 \\
 $V_{Ca}$ & 25 mV & $f_{er}$ & 0.01 \\
 ${\rm SERCA}_3$ & 0.2 ms$^{-1}$ & ${\rm SERCA_{2b}}$ & 0.02 $\mu$M ms$^{-1}$ \\
 $\beta$ & $4.5\times10^{-6}$ $\mu$M fA$^{-1}$ ms$^{-1}$ & $k_{PMCA}$ & 0.125 ms$^{-1}$ \\
 \hline
\end{tabular}
\caption{Parameter values used in all figures for the phantom bursting model, except $k_{PMCA}$, $r$, and $g_{K(Ca)}$ which are varied from their default values in Figs. \ref{fgr13: Phantom Burster gkca Low} and \ref{fgr14: Phantom Burster gkca High}.}
\label{table 2: Phantom burster parameters}
\end{table}

The dynamics of the membrane potential $V$ are governed by the different ionic currents and by the membrane capacitance $C_m$. The ionic current are given by:
\begin{equation*}
    \begin{split}
        I_{Ca} &= g_{Ca}m_{\infty}(V)(V-V_{Ca}), \\
        I_K &= g_Kw(V-V_K), \\
        I_{K(Ca)} &= g_{K(Ca)}\frac{c^5}{K_d^5+c^5}(V-V_K), \\
        I_{K(ATP)} &= g_{K(ATP)}a(V-V_K).
    \end{split}
\end{equation*}
The equilibrium activation functions are increasing Boltzmann functions:
\begin{equation*}
    \begin{split}
        m_{\infty}(V) = \frac{1}{1+\exp(\frac{v_m-V}{s_m})}, \\
        w_{\infty}(V) = \frac{1}{1+\exp(\frac{v_w-V}{s_w})}, \\
        a_{\infty}(c) = \frac{1}{1+\exp(\frac{r-c}{s_a})}. \\
    \end{split}
\end{equation*}
The last of these, for $a_\infty (c)$, contains the parameter $r$ that is varied in Fig. \ref{fgr13: Phantom Burster gkca Low}. Increasing $r$ shifts the sigmoidal curve to the right. 

The flux of Ca$^{2+}$ through the membrane is
\[
J_{mem} = -(\beta I_{Ca} + k_{PMCA}c),
\]
where $\beta$ converts units of current to flux, and $k_{PMCA}c$ is the flux through the plasma membrane Ca$^{2+}$ pumps. The parameter $k_{PMCA}$ is varied in Fig. \ref{fgr14: Phantom Burster gkca High}. We assume the Ca$^{2+}$ influx into the ER occurs via Sarcoplasmic Endoplasmic Reticulum Ca$^{2+}$ (SERCA) pumps with high (SERCA$_{2b}$) and lower affinity (SERCA$_3$):
\[
J_{\rm SERCA} = {\rm SERCA}_{2b} + {\rm SERCA}_3c.
\]
The efflux out of the ER is the Ca$^{2+}$ leak, which is taken to be proportional to the gradient difference between Ca$^{2+}$ in the cytosol and the ER is
\[
J_{leak} = p_{leak}(c_{er} - c).
\]
The net Ca$^{2+}$ efflux from the ER is 
\[
J_{er} = J_{leak} - J_{\rm SERCA}.
\]
The $c$ differential equation is then the sum of membrane and ER fluxes multiplied by the fraction of free total Ca$^{2+}$, $f_{cyt}$.

For the ER Ca$^{2+}$, $J_{er}$ is scaled by the ratio of the volumes of the cytoplasmic compartment ($V_{cyt}$) and the ER compartment ($V_{er}$) and multiplied by the fraction of free Ca$^{2+}$ in the ER ($f_{er}$) giving the differential equation for $c_{er}$.


\bibliographystyle{siam}
\bibliography{dynamic.bib}

\end{document}